# MULTIPLE SPACE DEBRIS COLLECTING MISSION OPTIMAL MISSION PLANNING

**Max CERF**

AIRBUS Defence and Space

78130 Les Mureaux, France

max.cerf@astrium.eads.net

## ABSTRACT

In order to keep a safe access to space in the coming years, it will be necessary to clean the LEO region from the most dangerous debris like spent satellites or launchers stages. An average removal rate of 5 debris per year is recommended to at least stabilize the current debris population. Successive missions must be planned over the years using similar vehicles in order to limit the development cost. This paper addresses the problem of the mission planning so that they can be achieved at minimal cost by a generic vehicle designed for such Space Debris Collecting missions.

The problem mixes combinatorial optimization to select and order the debris among a list of candidates, and continuous optimization to fix the rendezvous dates and to define the minimum fuel orbital maneuvers. The solution method proposed consists in three stages.

Firstly the orbital transfer problem is simplified by considering a generic transfer strategy suited either to a high thrust or a low thrust vehicle. A response surface modelling is built by solving the reduced problem for all pairs of debris and for discretized dates, and storing the results in cost matrices.

Secondly a simulated annealing algorithm is applied to find the optimal mission planning. The cost function is assessed by interpolation on the response surface based on the cost matrices. This allows the convergence of the simulated algorithm in a limited computation time, yielding an optimal mission planning.

Thirdly the successive missions are re-optimized in terms of transfer maneuvers and dates without changing the debris order. These continuous control problems yield a refined solution with the performance requirement for designing the future Space Debris Collecting vehicle.

The method is applicable for large list of debris and for various assumptions regarding the cleaning program (number of missions, number of debris per mission, total duration, deorbitation scenario, high or low thrust vehicle). It is exemplified on an application case with 3 missions to plan, each mission visiting 5 SSO debris to be selected in a list of 21 candidates.

# 1. Introduction

The near Earth region is crowded by space debris of all sizes. These debris originate from the old spacecrafts (satellites and launcher upper stages) released on orbit at the end of their operational life since the 1960[th]. The number of small debris grows constantly due to fragmentation or corrosion phenomena of these old spacecrafts. An efficient way to limit the proliferation is to remove the spent observation satellites mostly evolving on near-circular polar orbits in the altitude range 700-900 km altitude. Several studies recommend a removal rate of 5 heavy debris per year in order to stabilize the debris population[1,2,3,4].

A dedicated vehicle must be designed for such removal missions. This paper addresses the problem of planning the successive missions so that they can be achieved at minimal cost by similar vehicles.

## 1.1 Debris orbits

Most Low Earth Orbit (LEO) debris move on near circular orbits. At a given date $t_0$ a circular orbit is completely defined by its radius and two angles orientating the orbital plane in the Earth inertial reference frame. The classical orbital parameters are denoted (Figure 1) :

- $a(t_0)$ = semi-major axis (m)
- $I(t_0)$ = inclination (deg)
- $\Omega(t_0)$ = right ascension of the ascending node (RAAN) (deg)

The inclination I is the angle of the orbital plane with the Earth equatorial plane. The intersection of the orbital plane with the Equator is the line of nodes. The RAAN $\Omega$ is the angle between the X axis of the Earth inertial reference frame and the direction of the ascending node (node crossed with a northwards motion).

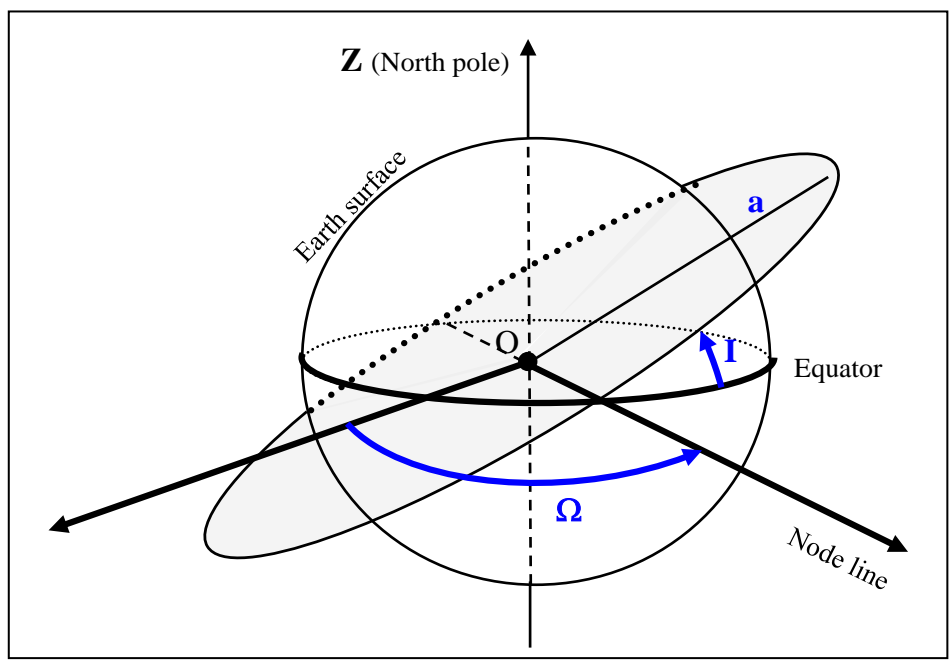

Figure 1 : Orbital parameters

The orbital parameters are constant in the keplerian model. The main perturbation to this model comes from the Earth flattening. Indeed the Earth equatorial bulge adds a perturbing force ($J_2$ zonal term) on the motion (represented by arrows on Figure 2). The resulting torque causes a precession of the orbital plane as pictured on the Figure 2.

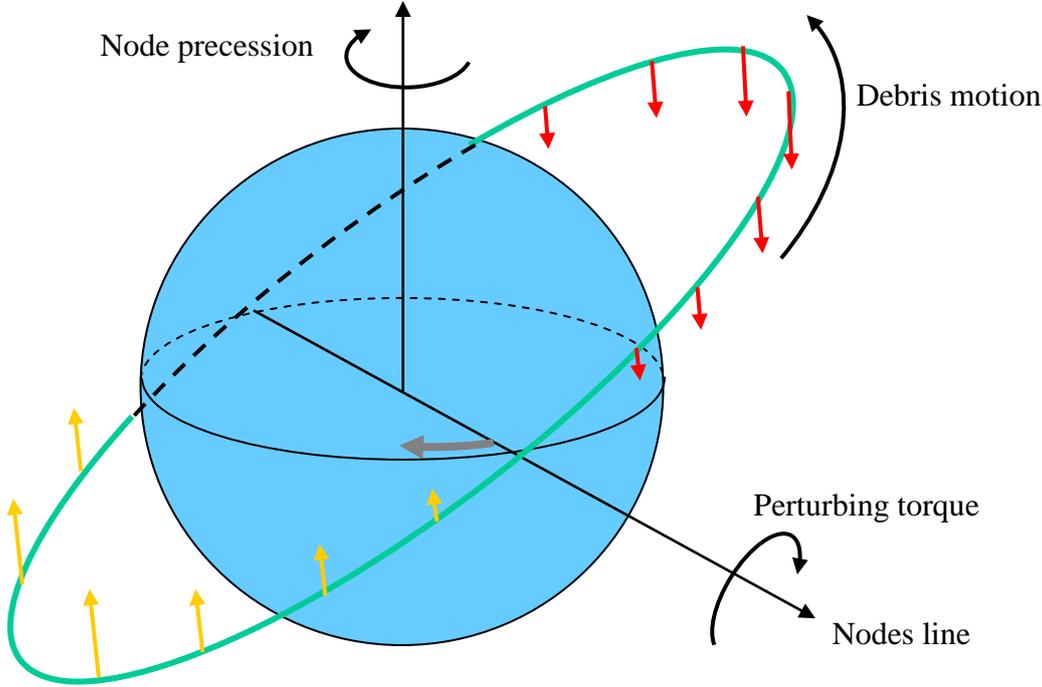

Figure 2 : Nodal precession due to the Earth flattening

The RAAN precession rate[9] depends on the orbit radius a and the inclination I :

$$\dot{\Omega} = -\frac{3}{2} J_2 \sqrt{\mu} R_T^2 \, a^{-\frac{7}{2}} \cos I$$

The constant of the Earth gravitational model are :

- $R_T$ = 6378137 m     (Earth equatorial radius)
- $\mu$ = 3.986005.$10^{14}$ m³/s²     (Earth gravitational constant)
- $J_2$ = 1.08266     (first zonal term)

The $J_2$ perturbation causes no secular change on the semi-major axis, the eccentricity and the inclination. The orbit remains circular with radius a and inclination I. The precession rate $\dot{\Omega}$ is therefore constant and the RAAN evolves linearly with the time :

$$\Omega(t) = \Omega(t_0) + \dot{\Omega}.(t - t_0)$$

For a sun-synchronous orbit (SSO), the RAAN precession rate matches the motion of the Sun direction as pictured on the Figure 3, with a rate equal to 0.986 deg/day (360 deg in 365.25 days). This property is favorable for an observation satellite since a region of given latitude is always flown over at the same local solar time. Most debris stemming from these spent observation satellites are consequently on nearly sun-synchronous orbits.

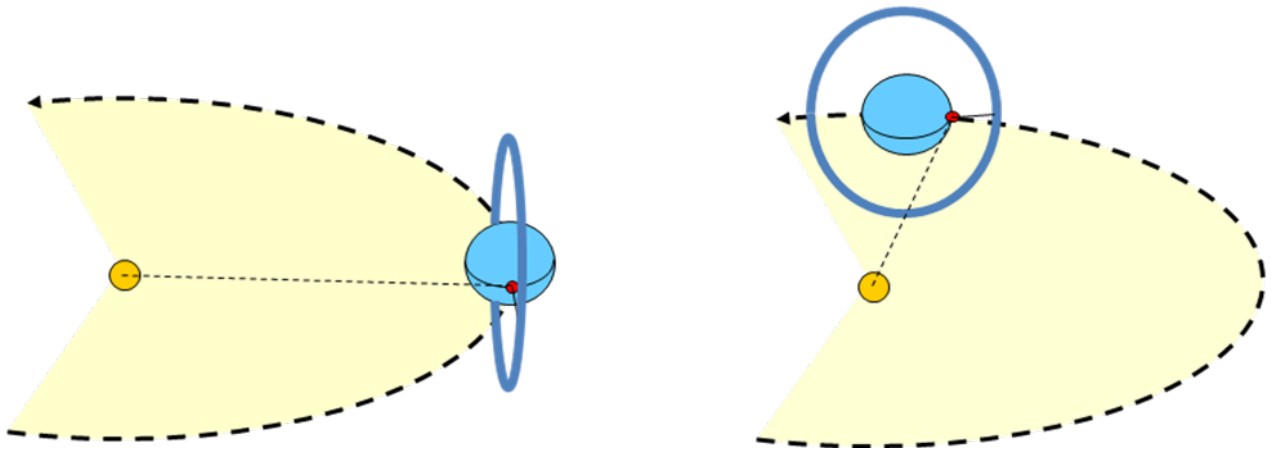

Figure 3 : Sun-synchronous orbit

## 1.2 Cleaning program

Several studies recommend a removal rate of 5 heavy debris per year from the LEO region of in order to stabilize the current population[1,2,3,4].

A cleaning program has to be defined to meet this requirement. It consists in launching a series of dedicated vehicles, each one being in charge of removing several debris (typically 5 debris per vehicle). With the current state of the art, a reusable concept cannot be envisioned. A series of expendable vehicles is necessary in order to progressively clean the LEO region from the most dangerous debris. In order to limit the development cost, the vehicles used for the successive missions should be similar.

The profile of a single Space Debris Collecting (SDC) mission is defined as follows :

- Choose the debris to visit
- Launch of the SDC vehicle onto orbit
- Travel from one debris to another
- Process each debris visited (observation, capture, deorbitation)
- Deorbit the vehicle itself at the end of the mission

The mission duration (typically 1 year) includes the transfers between the successive debris and the operations applied to each of them. The mission cost is driven at the first order by the SDC vehicle initial mass. This gross mass comprises the fuel required by the powered maneuvers (orbital transfers to go from one debris to another, deorbitation if performed by the vehicle itself) and the masses of the sub-systems used for the debris processing (rendezvous, capture, deorbitation if performed by an autonomous kit supplied to the debris).

The design of the SDC vehicle is based on the most expensive mission : this ensures that the same design is compliant of all the successive missions planned. The overall cost of the cleaning program is at the first order driven by the SDC vehicle gross mass.

## 1.3 Problem statement

The goal is to design a minimal mass vehicle compliant of a series of successive removal missions. The successive missions must be planned so as to minimize the fuel requirement of the most expensive mission, while achieving a mean removal rate of 5 debris per year. For that purpose, the following issues must be addressed :

- How to minimize the cost of a single mission (recurrent cost minimization) ?
- How to plan the successive removal missions (development cost minimization) ?

We denote :

- N the total number of debris considered in the list
- n the number of debris to visit per mission
- m the number of missions planned

The total number of debris visited at the end of the last mission is m×n out of N candidates. The m×n selected debris are visited at the successive dates $t_1$, $t_2$, ... , $t_{m \times n}$. The debris order and the rendezvous dates have to be optimized.

The SDC problem formulate as a graph problem. In terms of graph optimization, the debris are the nodes, the transfer trajectories are the edges while the successive missions are represented as opened sub-paths. The Figure 4 illustrates a 21 debris case, with 3 missions visiting each one 5 debris. Only 15 debris out of the 21 candidates will be visited, whilst 6 debris will be left on orbit.

- The 21 candidates debris are figured by the black points (nodes).
- The 3 missions are represented by respectively blue, green and red arrows (edges).
- Each mission deals with 5 debris (sub-paths).

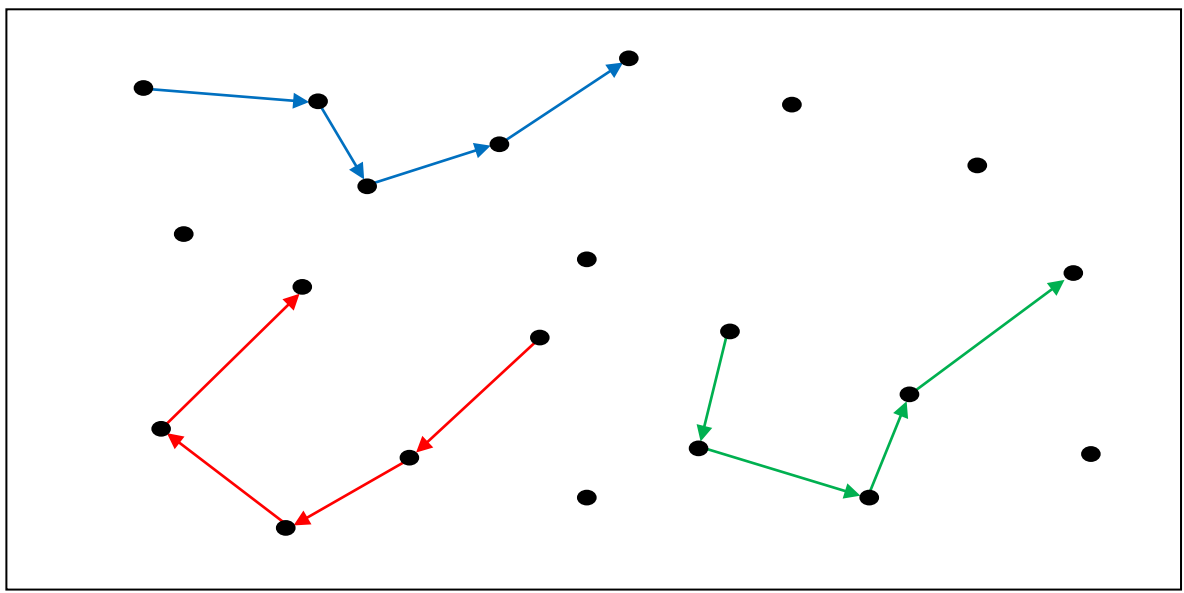

Figure 4 : Illustration of the SDC problem

The SDC problem is a variant of the Travelling Salesman Problem (TSP). The TSP consists in finding the minimal distance closed path visiting all the nodes once. The classical TSP features are the following :

- The nodes are fixed in a plane and the cost of going from one node to another is measured by the Cartesian distance in the plane (represented on the Figure 4 by the arrow length). The TSP problem is not time-dependent.
- Every node has to be visited once and once only and the path is closed. The overall cost is the path length.

There are three main differences between the TSP and the SDC :

- The Earth flattening causes the precession of the debris orbits plane. The precession rate is different between the debris, so that their relative configuration evolves with the time. The cost of going from the debris j to the debris k depends on the starting date $t_j$ and the arrival date $t_k$, making thus the SDC problem time-dependent.
- Instead of a single closed path visiting all the nodes, the debris are gathered in several sub-paths (missions). Not all debris are visited and the cost is measured from the most expensive sub-path.
- There is a global time constraint due to the targeted removal rate of 5 debris per year.

The cost evaluation procedure is illustrated on the Figure 5 in the case of 3 successive missions of 5 debris to be selected in a list of 21 candidates. The 15 selected debris are visited at the respective increasing dates $t_1 < t_2 < \ldots < t_{15}$. The respective costs of the 3 missions are $K_1$, $K_2$, and $K_3$. The cost K of the cleaning program is the cost of the most expensive mission. The total duration of the cleaning program $t_{15} - t_1$ must be lower than 3 years to achieve the targeted removal rate (5 debris per year).

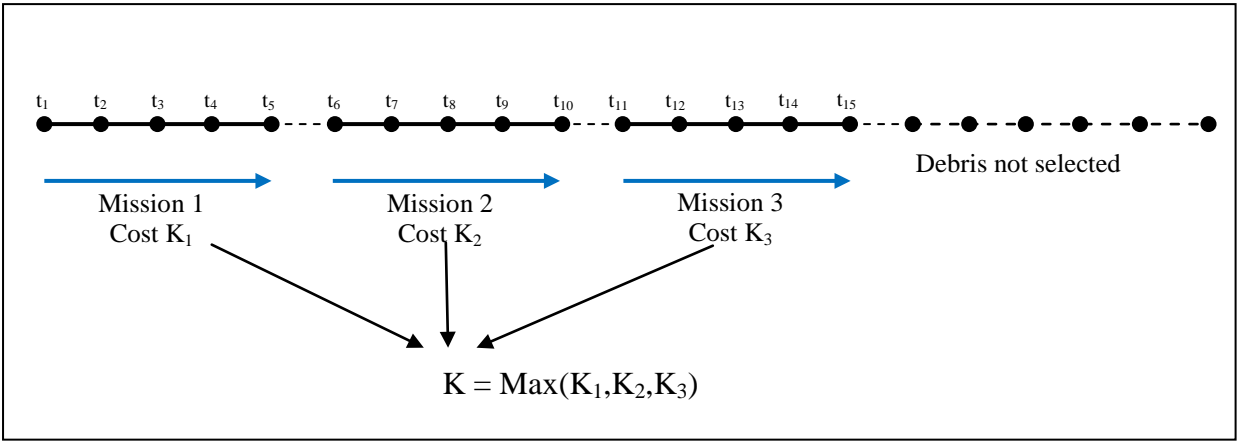

Figure 5 : SDC cost function

The differences between the TSP and the SDC are summarized in the Table 1.

|  | TSP | SDC |
|---|---|---|
| Number of nodes visited | N | m×n ≤ N |
| Path definition | Single closed path | Several opened sub-paths |
| Node positions | Fixed | Moving |
| Edge valuations | Fixed length | Time-dependent |
| Cost function | Path length | Maximum sub-path cost |
| Time constraint | None | Upper bound |

Table 1 : TSP vs SDC

A major issue in the SDC problem lies in the valuation of the edges. Each edge represents the orbital transfer between a debris and the next one on the path. The edge valuation is the propellant required to perform the orbital transfer. Finding this minimal fuel trajectory is a challenging optimal control problem.

The global problem consists thus in a series of continuous problems (transfer trajectories between debris) embedded within a combinatorial problem (path between the selected debris). It mixes integer variables (debris selection and order), real variables (rendezvous dates) and optimal control (transfer maneuvers). Even taken separately these sub-problems are intrinsically hard. It is out of reach to solve the global SDC problem in a direct manner.

### 1.4  Optimization method

The global problem is a complex variant of the TSP which is a NP-complete problem. Instances of the TSP are used as benchmark for combinatorial optimization algorithms. A large number of algorithms have been experimented[5,6,7]. They can be roughly classified as presented in Table 2.

| Algorithms | Explicit enumeration | Implicit enumeration | Greedy heuristics Stochastic programming |
|---|---|---|---|
| Solution accuracy | Exact | Exact | Approximate |
| Issues | CPU time | Linearization + Iteration | Algorithm settings |
| Problem size | Small | Medium | Large |

Table 2 : Algorithms for combinatorial optimization

Finding the exact solution of a combinatorial problem requires an enumeration algorithm, either explicit or implicit. Such algorithms are applicable only to limited size problems[8]. For large instances, only approximate solutions can be hoped in a fixed computation time. In view of handling large lists of candidate debris, we have to turn to stochastic algorithms[7]. Among all the existing algorithms, simulated annealing has proved quite successful on large TSP instances. We have therefore selected a simulated annealing approach to tackle the SDC problem.

Compared to the TSP, the SDC problem presents additional issues due to the edge valuations and their time-dependency. Indeed a simulated annealing algorithm tries millions of solutions before achieving a satisfactory convergence. Each trial solution is defined by a debris order and the visiting dates. Assessing the exact cost function (measured by the fuel consumption) of a trial solution requires solving a series of hard optimal control problems to find the trajectories between the successive debris.

In order to apply a simulated annealing to the SDC problem with reasonable computation times, it is not possible to solve "on-line" these optimal control problems for each trial solution. An instantaneous cost function must be devised. In order to get a sufficient confidence in the simulated annealing results, this cost function must be both robust (i.e. yield a cost value whatever the input data) and reliable (i.e. yield a cost value representative of a real optimized transfer).

The approach proposed consists in using a Response Surface Modelling (RSM) based on cost matrices. More precisely the optimization process is split into three successive stages.

- The first stage consists in building the cost matrices. These cost matrices store the costs of all the possible elementary transfers between debris for a mesh of discretized dates. They result from a series of optimizations based on a simplified generic transfer strategy adapted to the mission specificities and to the vehicle propulsion system (high thrust or low thrust).
- The second stage consists in finding the optimal mission planning with a simulated annealing algorithm. The algorithm is derived from the one applied to a classical TSP with additional variables (rendezvous dates) and a RSM based cost function. The cost of a trial solution is assessed by interpolation in the cost matrices spanning the possible transfers and dates. The simulated annealing solution defines the optimal mission planning.
- The third stage consists in a refined trajectory optimization. Indeed the RSM yields an approximate cost value by interpolation in the cost matrices. The refined optimization consists in fixing the debris order for the successive missions as given by the simulated annealing, and optimizing the rendezvous dates and the maneuvers using a real trajectory simulation. This is a standard optimal control problem with continuous variables. The solution yields the requirement (fuel or velocity impulse) for the vehicle design.

The optimization process is pictured on the Figure 6 with the algorithms used for the three optimization stages.

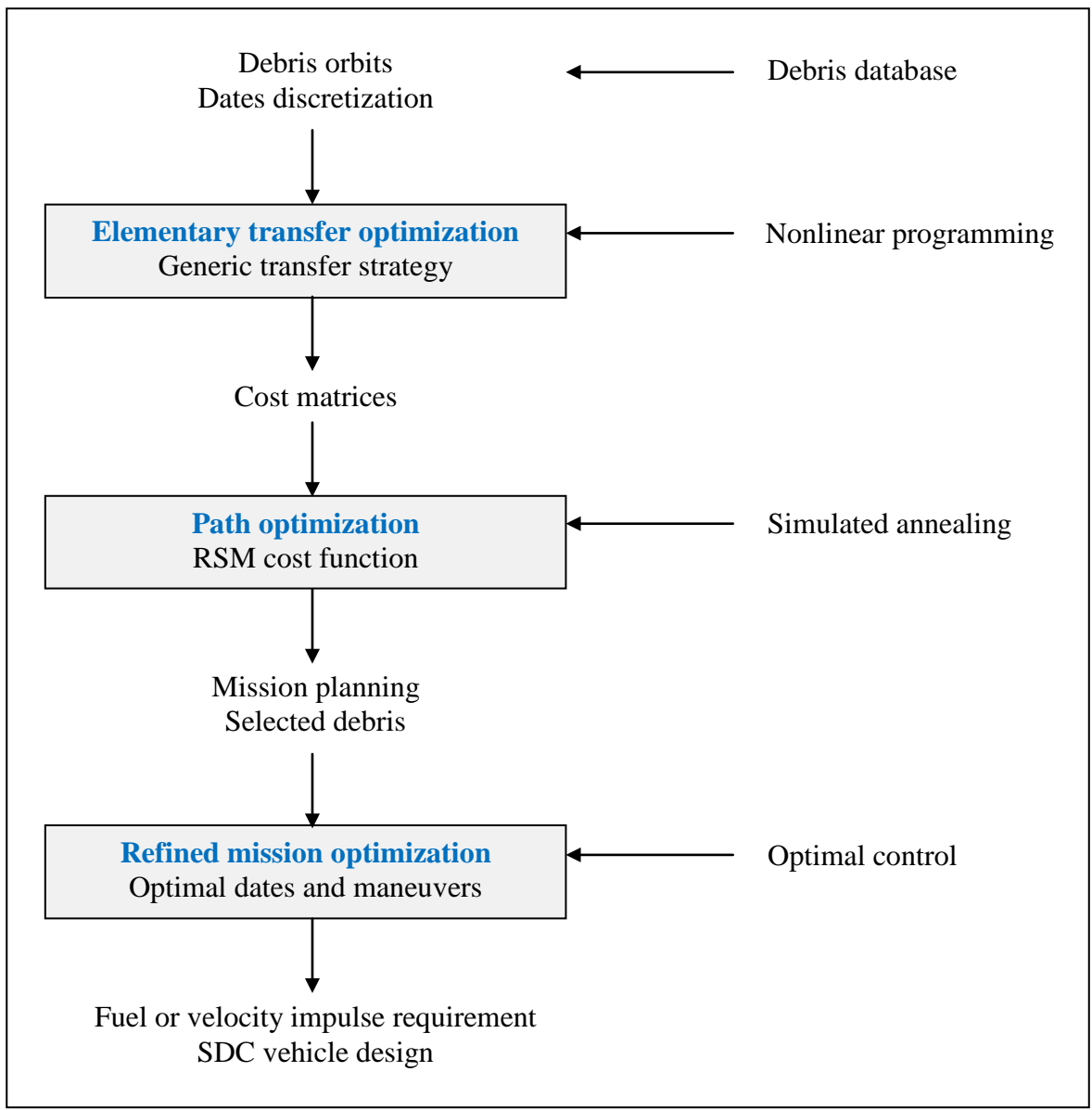

Figure 6 : Optimization process

The next sections detail the solution methods proposed for the transfer problem (§2) and the path problem (§3). The practical implementation of the overall process is presented (§3) and illustrated on an application case (§4).

## 2. Transfer problem

Finding the minimal fuel trajectory from a debris to another is a difficult optimal control problem in the general case. This transfer problem is by considering a generic transfer strategy adapted to the mission specificities and to the vehicle propulsion system. The optimal control problem reduces thus to a nonlinear programming problem with two variables and one constraint that can be solved in an efficient manner. This simplified modelling is used to build the cost matrices used in the RSM cost function.

### 2.1 Transfer strategy

The simplifications of the transfer problem are based on the mission specificities :

- The orbits of the targeted debris (old observation satellites) are assumed to be circular. The real orbits of such debris have indeed negligible eccentricities (e < 0.01).
- The mean removal rate (5 debris per year) allocates an average duration of 3 months per transfer. This duration leaves time enough to use the $J_2$ nodal precession in order to perform the RAAN change at null fuel consumption.

The generic transfer strategy consists in bringing the vehicle on a circular drift orbit and wait until the RAAN change is completed. More precisely the transfer from a debris 1 to a debris 2 is split into three phases :

- A propelled transfer from the debris 1 orbit to the drift orbit
- A waiting duration on the drift orbit
- A propelled transfer from the drift orbit to the debris 2 orbit

The transfer starts at a given date $t_1$ and ends at a given date $t_2$. The orbital parameters of the successive orbits are denoted in the Table 3.

|  | Propelled transfer 1 | | Drift phase | | Propelled transfer 2 |
|---|---|---|---|---|---|
| Orbit | Debris 1 | Drift start | | Drift finish | Debris 2 |
| Date | $t_1$ | $t_{d1}$ | | $t_{d2}$ | $t_2$ |
| Radius | $a_1$ | $a_d$ | | $a_d$ | $a_2$ |
| Inclination | $I_1$ | $I_d$ | | $I_d$ | $I_2$ |
| RAAN debris 1 | $\Omega_1(t_1)$ | $\Omega_1(t_{d1})$ | | $\Omega_1(t_{d2})$ | $\Omega_1(t_2)$ |
| RAAN debris 2 | $\Omega_2(t_1)$ | $\Omega_2(t_{d1})$ | | $\Omega_2(t_{d2})$ | $\Omega_2(t_2)$ |
| RAAN vehicle | $\Omega_v(t_1) = \Omega_1(t_1)$ | $\Omega_v(t_{d1})$ | | $\Omega_v(t_{d2})$ | $\Omega_v(t_2) = \Omega_2(t_2)$ |

Table 3 : Successive orbits during the transfer

The rendezvous in anomaly with the debris 2 is neglected both in terms of duration and consumption compared to the overall transfer. This generic transfer strategy using the $J_2$ precession to control the RAAN at null fuel consumption is near-optimal as long as a sufficient duration ($t_2 - t_1$) is allocated. For short durations this strategy would no longer be possible and the RAAN change should be realized by propelled maneuvers at the expense of a larger fuel consumption. Two modellings of the propelled transfers are considered depending whether the SDC vehicle uses a high thrust or a low thrust propulsion system.

### 2.1.1 High thrust propulsion

In the case of a high thrust engine the powered orbital transfers are modelled as Hohmann transfers with impulsive maneuvers. Each orbital transfer (from debris 1 to drift, then from drift to debris 2) is achieved by a two impulse Hohmann transfer[9,10] with split inclination change. The inclination of the intermediate elliptical orbit is computed using a near-optimal approximation derived by Lisowski[10]. The approximation consists in minimizing the sum of the squared velocity impulses (instead of the velocity impulses norm). An analytical solution can thus be found with a limited deviation from the true minimum.

The transfer strategy is depicted on the Figure 7, with the successive velocity impulses associated to the initial and final Hohmann transfers.

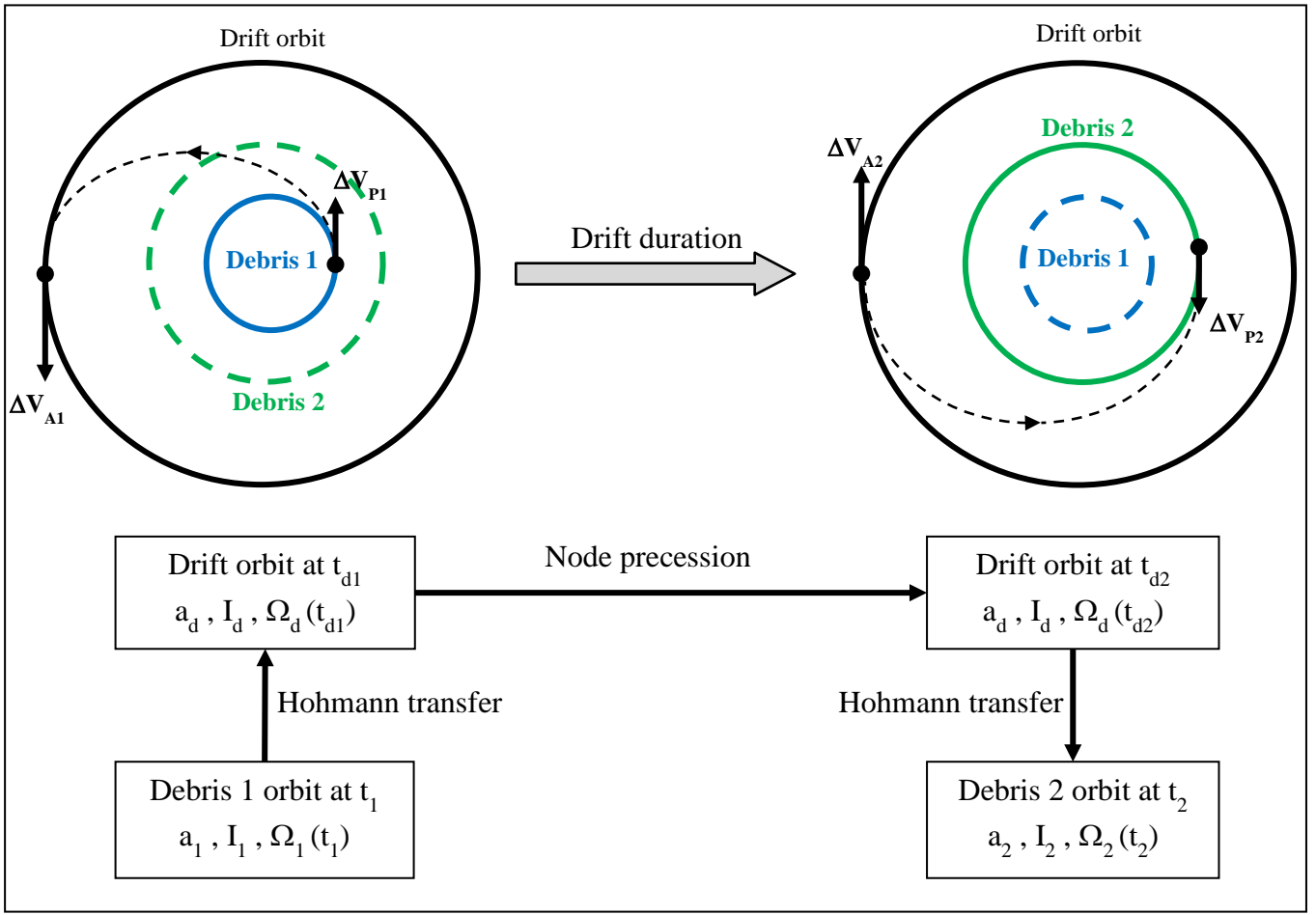

Figure 7 : High thrust transfer strategy

The Hohmann transfer durations (about 1h) are negligible wrt the drift duration (several days or weeks). The RAAN precession due to the $J_2$ may be neglected during these transfers :

$$\begin{cases} t_{d1} \approx t_1 & \Rightarrow \quad \Omega_v(t_{d1}) \approx \Omega_v(t_1) \\ t_2 \approx t_{d2} & \Rightarrow \quad \Omega_v(t_2) \approx \Omega_v(t_{d2}) \end{cases}$$

The transfer total cost is measured by summing the four velocity impulses. It does not depend on the vehicle thrust level.

### 2.1.2 Low thrust propulsion

In the case of a low thrust engine the powered orbital transfers are modelled as Edelbaum transfers with continuous thrusting. Each orbital transfer (from debris 1 to drift, then from drift to debris 2) is achieved by a minimum time Edelbaum transfer with continuous inclination change[11,12]. The Edelbaum model assumes a constant acceleration level. In order to get a refined assessment of the transfer duration and cost, the solution is computed in two stages : first with the initial acceleration level, then with the average acceleration level estimated from the first solution.

The transfer strategy is depicted on the Figure 8, with the spiraling trajectories associated to the initial and final Edelbaum transfers.

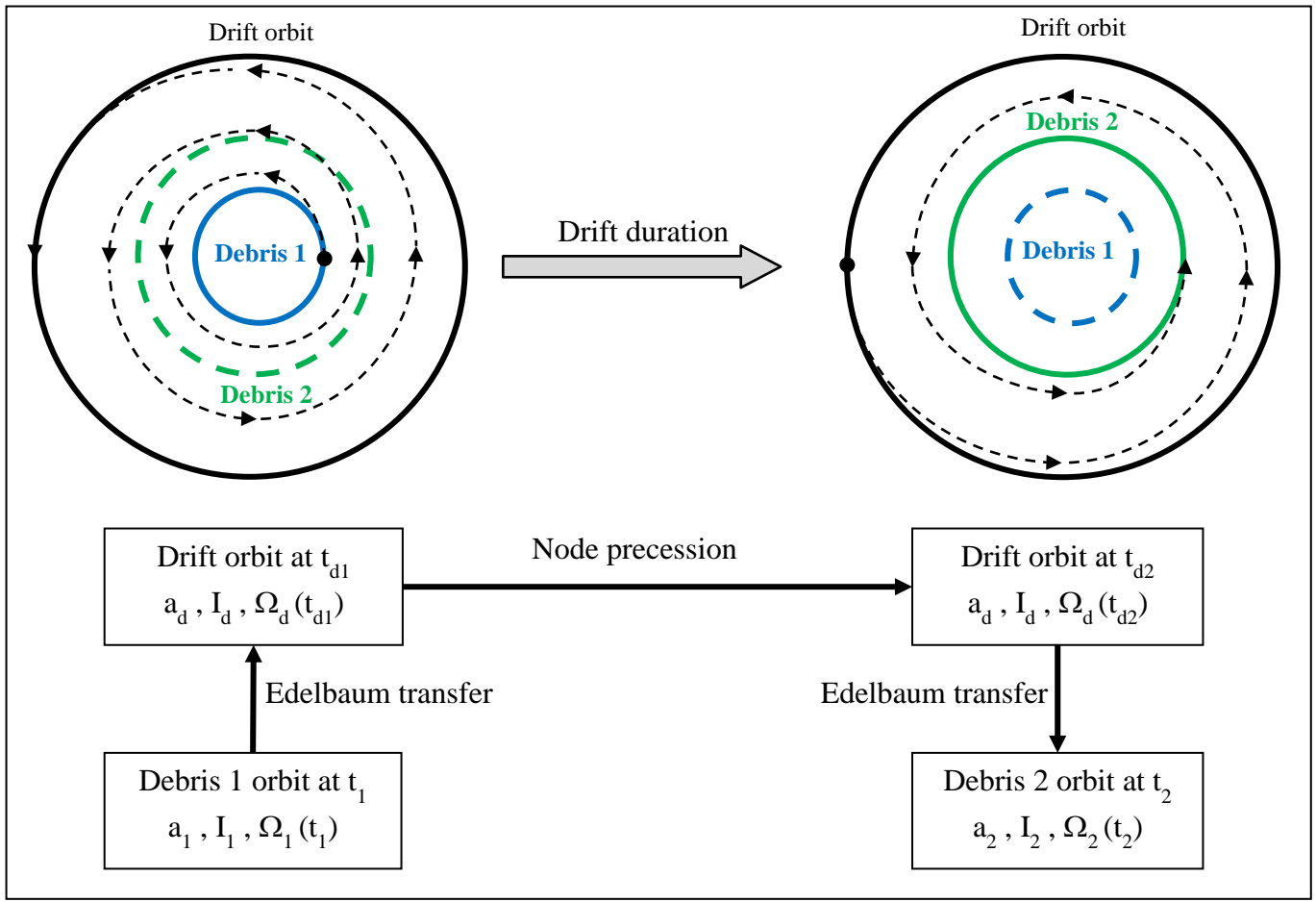

Figure 8 : Low thrust transfer strategy

Opposite to the high thrust case, the durations of the Edelbaum transfers are no longer negligible wrt the drift duration and they may induce significant RAAN changes. The RAAN evolution is assessed by a numerical integration along the Edelbaum trajectory :

$$\begin{cases} \Omega_v(t_{d1}) = \Omega_v(t_1) + \int_{t_1}^{t_{1d}} \dot{\Omega}_v(t)dt & \text{during the first propelled transfer} \\ \Omega_v(t_2) = \Omega_v(t_{d2}) + \int_{t_{d2}}^{t_2} \dot{\Omega}_v(t)dt & \text{during the second propelled transfer} \end{cases}$$

The Edelbaum solution yields the minimal time transfer between mutually inclined circular orbits, assuming a constant acceleration level. The Edelbaum model is based on an averaging of the dynamic equations assuming that the orbit remains circular throughout the transfer. During each revolution, the thrust direction keeps a constant angle with the orbital plane, with a sign switch at the antinodes. This averaged control law does not modify directly the RAAN. The RAAN evolution is only due to the $J_2$ perturbation which acts constantly throughout the transfer phases.

The Edelbaum solution yields the evolution of the mean orbit radius a(t) and inclination I(t) throughout the transfer. The mean RAAN precession rate is computed as :

$$\dot{\Omega}(t) = -\frac{3}{2} J_2 \sqrt{\mu} R_T^2 \, a(t)^{-\frac{7}{2}} \cos I(t)$$

The RAAN variation during the propelled transfer is assessed by a numerical integration from the initial date $t_1$ to the final date $t_2$.

The velocity impulse associated to the Edelbaum solution is obtained as the product of the mean acceleration level denoted f by the transfer duration $t_2 - t_1$ :   $\Delta V = f.(t_2 - t_1)$.
The transfer total cost is measured by summing the velocity impulses of the two propelled transfers (from debris 1 to drift, then from drift to debris 2). Opposite to the high thrust case, this cost depends on the vehicle thrust level.

Remark

The transfer strategy based on Edelbaum transfers is not globally optimal. Indeed the Edelbaum solution yields the minimal time transfer without taking into account the RAAN change. The RAAN change is assessed a posteriori along the Edelbaum trajectory. The drift orbit parameters $a_d$ and $I_d$ must then yield the adequate precession $\dot{\Omega}_d$ rate to achieve the required RAAN final value.

This may lead to a more costly drift orbit regarding the velocity impulses. Cheaper solutions could be found by performing a part of the RAAN change during the propelled transfers. The possible cost gain may be significant depending on the relative durations of the propelled phases wrt the drift phase. Variants of the Edelbaum solution have been derived considering alternative constraints[12] (RAAN change instead of inclination change, altitude bound). For the SDC problem an analytical solution taking into account the three transfer phases (propelled – drift – propelled) is currently under investigation.

It is assumed for the SDC problem that a sufficient acceleration level is available on the vehicle and that a sufficient transfer duration is allocated so that the transfer strategy using the $J_2$ is near optimal. With these assumptions the propelled durations should remain small wrt the drift duration, and the Edelbaum transfer strategy can be considered as nearly optimal.

## 2.2 Problem formulation

The transfer optimization consists in finding the drift orbit semi major axis $a_d$ and inclination $I_d$ in order to :

- Attain the RAAN of the debris 2 → constraint
- Minimize the fuel consumption → cost function

The transfer starts at the date $t_1$ and finishes at the date $t_2$. These initial and final dates are fixed.

### 2.2.1 Constraint

The vehicle and debris 2 RAAN are denoted respectively $\Omega_v(t)$ and $\Omega_2(t)$. Their evolution is only due to the $J_2$ perturbation which acts constantly throughout the transfer phases.

The vehicle RAAN must go from the debris 1 RAAN $\Omega_1(t_1)$ at the transfer beginning to the debris 2 RAAN $\Omega_2(t_2)$ at the transfer end. The RAAN constraint is thus expressed as :

$$\Omega_v(t_2) = \Omega_2(t_2) \quad \text{with} \quad \begin{cases} \Omega_v(t_2) = \Omega_1(t_1) + \int_{t_1}^{t_{1d}} \dot{\Omega}_v(t)dt + \dot{\Omega}_d \cdot (t_{d2} - t_{d1}) + \int_{t_{d2}}^{t_2} \dot{\Omega}_v(t)dt \\ \Omega_2(t_2) = \Omega_2(t_1) + \dot{\Omega}_2 \cdot (t_2 - t_1) \end{cases}$$

### 2.2.2 Cost function

The fuel consumption comes from the propelled transfer respectively from the debris 1 to the drift orbit, and from the drift orbit to the debris 2.

The mass consumed $m_c$ is linked to the velocity impulse by the rocket equation[9,10] :

$$\Delta V_{12} = v_e \ln \frac{M_1}{M_2} \quad \Leftrightarrow \quad m_c = M_1 - M_2 = M_1 \left( 1 - e^{-\frac{\Delta V_{12}}{v_e}} \right)$$

where
- $v_e$ is the exhaust velocity of the vehicle engine
- $M_1$ is the vehicle mass at the transfer beginning (date $t_1$)
- $M_2$ is the vehicle mass at the transfer end (date $t_2$)

For a given initial mass $M_1$, minimizing the fuel consumption is equivalent to minimizing the velocity impulse. The velocity impulse is preferred as cost function rather than the mass consumed. It is indeed an intrinsic cost measure independent on the transfers previously realized by the vehicle, opposite to the mass consumption which depends on the vehicle gross mass at the transfer beginning. It will therefore be more suited to the optimization of successive transfers required by the SDC mission.

### 2.2.3 NLP problem

The transfer optimization problem is formulated as:

$$\min_{a_d, i_d} \Delta V \text{ s.t. } \Omega_v(t_2) = \Omega_2(t_2)$$

The optimization variables are the drift orbit radius $a_d$ and inclination $I_d$. The initial and final dates, respectively $t_1$ and $t_2$, are fixed.

In some cases it can be more economical to complete the transfer at a date prior to $t_2$, for example when the precession rate on the initial orbit is sufficient to naturally compensate the RAAN difference between the debris 1 and the debris 2 within the allocated duration. In such cases, the transfer on an intermediate drift orbit is useless. The vehicle waits on the initial orbit (debris 1) and the transfer to the debris 2 orbit takes place when the RAAN difference is nullified. The final date is then lower than $t_2$. These cases are accounted in the above formulation by allowing a drift orbit identical to the debris 1 orbit, and by adding an optional waiting phase on the debris 2 orbit until reaching the fixed final date $t_2$.

This reduced optimization problem with 2 variables and 1 constraint is readily solved with a nonlinear optimizer, taking into account either the Hohmann transfer strategy for a high thrust engine, or the Edelbaum transfer strategy for a low thrust engine. Since the debris RAAN change with the time, the minimal cost denoted $\Delta V_{opt}$ is a nonlinear function of the starting and arrival dates, respectively $t_1$ and $t_2$, or equivalently of the starting date $t_1$ and the transfer duration $\Delta t = t_2 - t_1$.

Remark

Some care must be taken regarding the initialization of the optimization variables. Indeed the nonlinear problem has at least two local minima depending on whether the debris 2 RAAN is attained forwards or backwards. The drift precession rate $\dot{\Omega}_d$, which depends on the optimization variables $a_d$ and $I_d$, must be initialized correctly in order to converge on the best of these two minima. The adequate rate is chosen depending on the respective RAAN of the debris 1 and 2 at the transfer beginning, as pictured on the Figure 9.

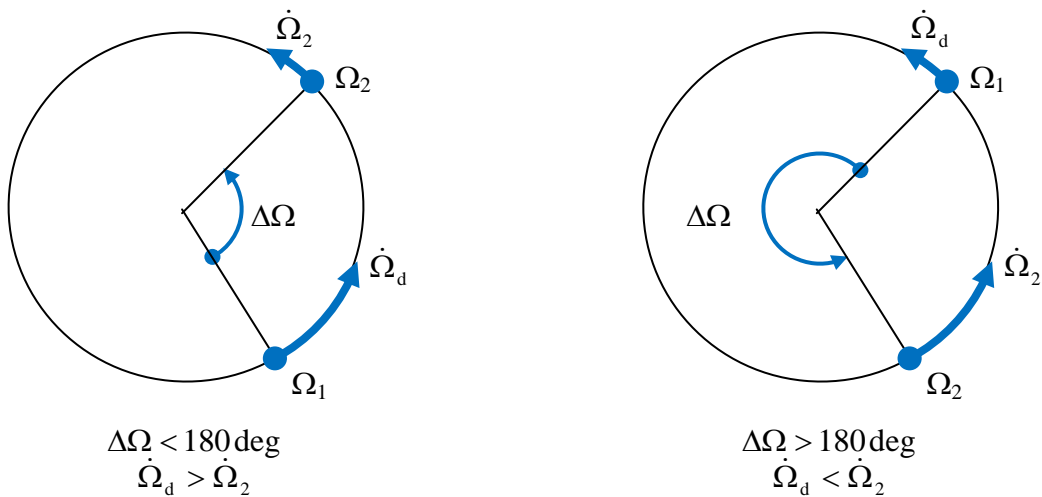

Figure 9 : Forwards (left) or backwards (right) correction

## 2.3 Cost matrices

Although the reduced problem can be solved in a robust and efficient manner by a nonlinear optimizer, the computation still requires a few seconds. An on-line optimization is not suited to a simulated annealing algorithm which needs millions of trial to converge on an acceptable solution. In order to get reasonable computation times, the on-line optimization is replaced by a response surface approach based on cost matrices.

The cost of an elementary transfer from any debris to any other depends on the starting date and the transfer duration. For a given starting date $\tau$ and a given transfer duration $\Delta\tau$, it is possible to assess the transfer costs between all the pairs of the N candidate debris. This requires $N \times (N-1)$ elementary optimizations whose results are stored in a $N \times N$ cost matrix represented on the Figure 10.

| Date $\tau$ <br> Duration $\Delta\tau$ |  | Debris k <br> ($1 \leq k \leq N$) |  |
|---|---|---|---|
|  |  |  |  |
| Debris j <br> ($1 \leq j \leq N$) |  | Cost $C(\tau,\Delta\tau,j,k)$ |  |
|  |  |  |  |

Figure 10 : Cost matrix for the date $\tau$ and the duration $\Delta\tau$

The value $C(\tau,\Delta\tau,j,k)$ stored at the row j and column $k \neq j$ of the matrix is the optimal cost $\Delta V_{opt}$ to go from the debris j at the date $\tau$ to the debris k at the date $\tau+\Delta\tau$. The matrix diagonal is unfilled at this stage. It will be used later (§3) to account for the cost of the debris operations.

In order to account for the time-dependency of the SDC problem, a series of cost matrices are assessed for a mesh of discretized starting dates and transfer durations covering the time span of the cleaning program. We denote :

- $T_0$  the date of the beginning of the cleaning program
- $\Delta T$ the total duration of the cleaning program
- $n_t$   the number of discretized starting dates
- $n_d$   the number of discretized transfer durations
- $\tau_i$   the starting date number i in the grid          ($1 \leq i \leq n_t$)
- $\Delta\tau_d$ the transfer duration number d in the grid     ($1 \leq d \leq n_d$)

For any starting date $\tau_i$, any duration $\Delta\tau_d$ and any pair of debris j and $k \neq j$, $C(i,d,j,k)$ is the cost of the transfer going from the debris j at the date $\tau_i$ to the debris k at the date $\tau_i+\Delta\tau_d$.
The sub-matrix $C(i,d,1:N,1:N)$ of size $N \times N$ contains the costs of all the elementary transfers starting at the date $\tau_i$ with a duration $\Delta\tau_d$. It requires $N \times (N-1)$ optimizations for solving the associated transfer problems. Some transfers may be unfeasible in the prescribed duration due to bounds on the drift orbit parameters (minimal altitude) that limit the available precession rate. In such cases, the corresponding matrix element is set to an arbitrarily large value, so that it will not be selected during the path optimization.

The total number of N×N sub-matrices is $n_t \times n_d$, corresponding to the mesh of $n_t$ dates and $n_d$ durations. Some of these sub-matrices are theoretically useless when they corresponds to a final date ($\tau_i + \Delta\tau_d$) beyond the ending date of the cleaning program ($T_0 + \Delta T$). These useless matrices are not computed and filled with large values indicating the transfer unfeasibility.

The Figure 11 illustrates the mesh of $n_t \times n_d$ cost sub-matrices spanning the dates of the cleaning program.

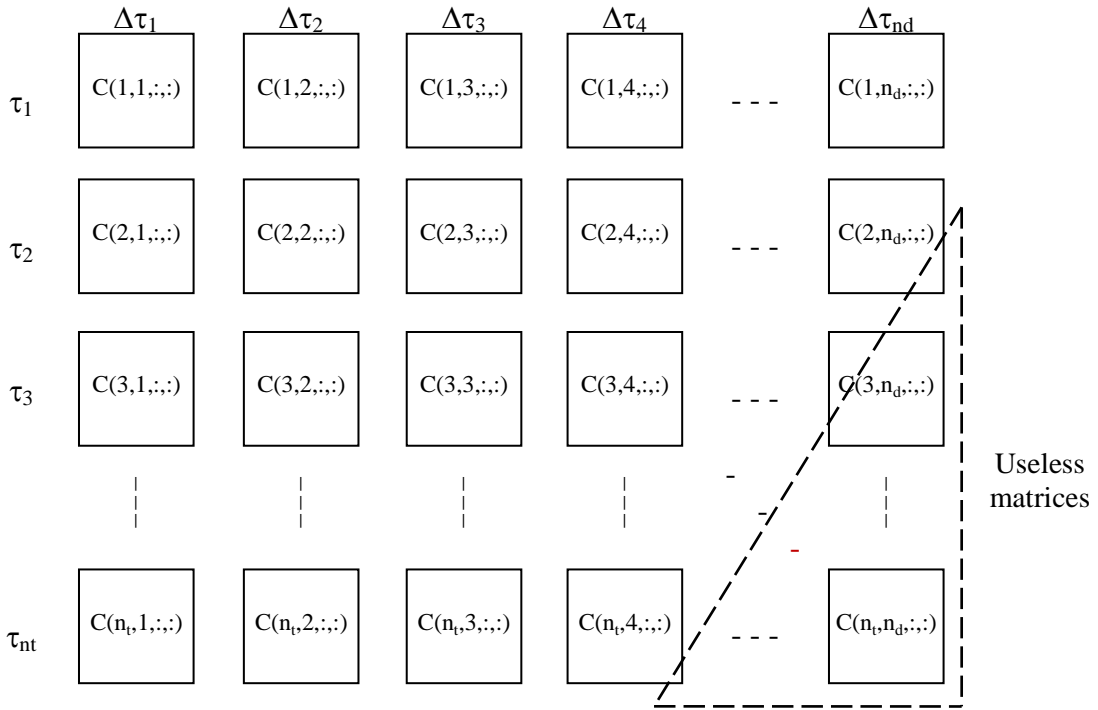

Figure 11 : Mesh of discretized cost matrices

A total of $n_t \times n_d \times N \times (N-1)$ optimizations is necessary to fill completely the mesh of cost matrices. This mesh is used within the simulated annealing process to assess the cost function through a Response Surface Modelling (RSM).

The RSM consists in a bilinear interpolation on the actual starting date and the actual duration. In order to avoid extrapolations that could lead to erroneous cost assessments, it is necessary to keep some "bounding matrices" in the mesh, particularly :

- A last row with a starting date greater than the ending date of the cleaning program. This lead to choose as last starting date : $\tau_{nt} = T_0 + \Delta T$
- On each row (with a starting date $\tau_i$) either the maximal transfer duration ($\Delta\tau_{nd}$), or the smallest transfer duration $\Delta\tau_d$ exceeding the ending date of the cleaning program ($\tau_i + \Delta\tau_d > T_0 + \Delta T$).

## 3. Path problem

The path problem consists in defining m successive missions (sub-paths) visiting each one n debris chosen among a list of N candidates. This graph problem is a time-dependent variant of the Travelling Salesman Problem (TSP). This section presents the solving method based on a simulated annealing algorithm, with a cost assessment method based on a response surface modelling.

### 3.1 Simulated annealing

Annealing is a metallurgic process to get an alloy without default. It consists in first melting the metal. At high energy level the atoms move freely and can exchange their positions. The metal is then cooled down very slowly. When their energy level decreases the atoms tend to freeze and to order in a crystalline structure. The quality of the alloy depends on the temperature decrease rate.

Simulated annealing is a stochastic optimization algorithm inspired by this metallurgic process[6,7]. It has been applied successfully to combinatorial problems with a large number of local minima, and particularly to the TSP. The algorithmic principles are the following :

- The current solution is noted $x_0$, its cost $f_0$ represents the energy level of the solution.
- A random perturbation is applied on $x_0$, yielding a candidate solution x with cost f.
- The candidate solution is accepted with the probability P computed as :

$$P = \begin{cases} 1 & \text{if } f \leq f_0 \\ e^{\frac{f_0 - f}{T}} & \text{if } f > f_0 \end{cases} = \text{probability of transition from the energy level } f_0 \text{ to } f$$

A degrading solution can therefore be accepted with a probability depending on the "temperature" parameter T : the higher the temperature, the higher the acceptance probability. The temperature threshold acts as an energy barrier which restricts the possibility of degrading the current solution (Figure 12).

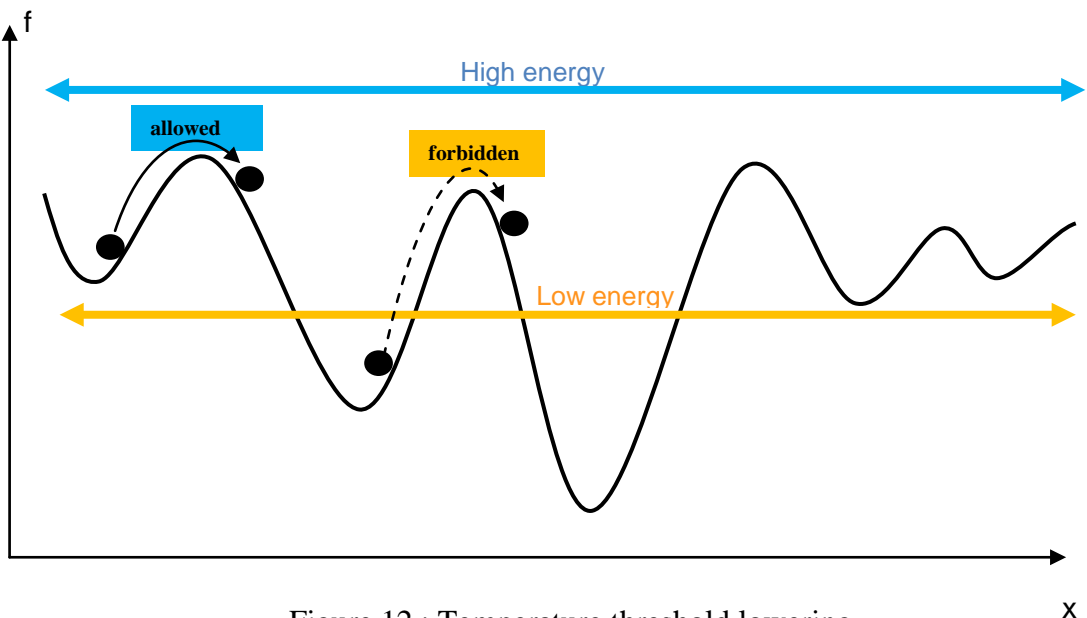

Figure 12 : Temperature threshold lowering

This mechanism allows escaping local minima by accepting random uphill moves and exploring widely the cost function landscape. When the temperature is progressively lowered, the solution freezes on the best minimum found.

The main settings of the algorithm are :

- The initial temperature $T_0$, the decrease rate $\alpha<1$ ($T_{k+1} = \alpha T_k$) and the number of tries at each temperature threshold.
- The definition of the random perturbations (or moves) applied to the current solution.

For each application case different values of the temperature parameters $T_0$ and $\alpha$ must be tried to get a satisfactory convergence. A too fast temperature decrease may trap the solution in a local minimum, while a too slow temperature decrease may result in a too large computation time.

Four elementary moves are implemented for the SDC problem : insertion, swap, permutation, date shift (Figure 13).
The insertion, swap and permutation modify the debris order on the path[7]. The date shift only changes the date of a node while keeping the path order. The new date remains comprised between the previous and next node date.

A single evaluation (or try) consists in :

- Selecting randomly one of the 3 elementary path moves (insertion, swap, permutation)
- Selecting randomly the nodes where the move is applied
- Performing the move to get the trial path
- Selecting randomly a node on the path
- Shifting randomly the node date between the previous and the next node date
- Assessing the cost of the trial solution
- Accepting the try with the probability level defined by the current temperature

An iteration of the simulated annealing consists in decreasing the temperature with a fixed rate $\alpha$ : $T_{k+1} = \alpha T_k$ after a given number of tries. A typical decrease rate value is $\alpha= 0.999$ every 1000 tries (this depends on the problem size).
The algorithm is initiated either with a random solution or with a greedy solution. For example, the initial solution can be built by the best insertion method : the nodes are inserted successively to the path at the position minimizing the cost. The initial temperature $T_0$ is set in order to accept a random perturbation of the initial solution with a 90% probability. This acceptation level is high enough to allow large solution changes during the first iterations.

When no progress is made after several temperature thresholds, a local search is performed by trying systematically all the elementary perturbations on the last solution. If this search is successful, the iterations are retrieved from the improved solution, else the algorithm is stopped. The overall algorithm is depicted on the Figure 14.

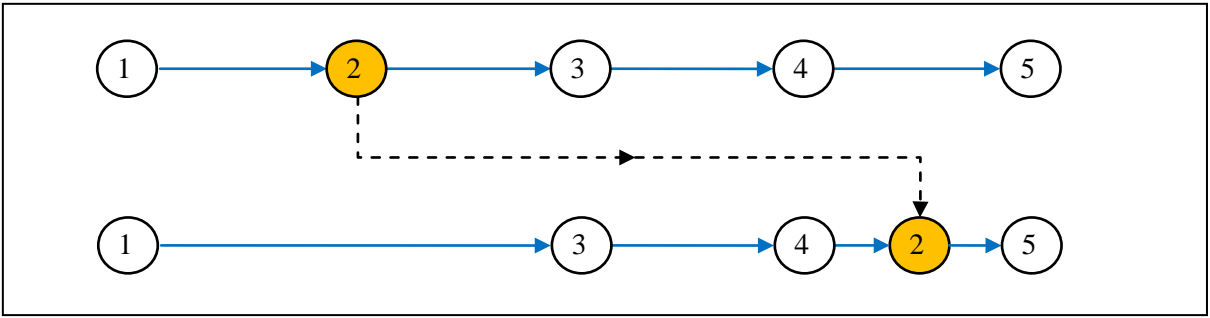

Figure 13.1 : Insertion (2 is inserted after 4)

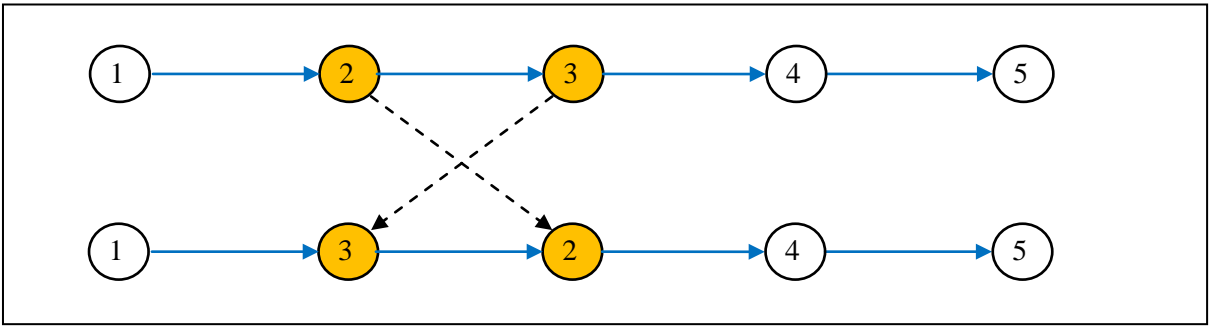

Figure 13.2 : Swap (2 and 3 are exchanged)

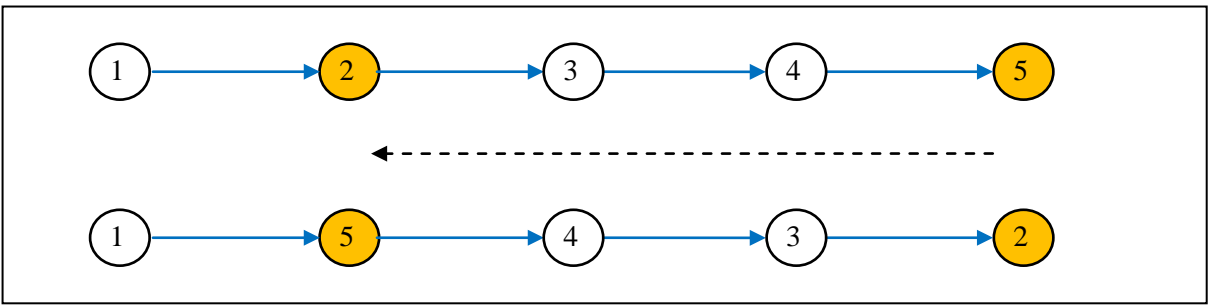

Figure 13.3 : Permutation (the leg from 2 to 5 is reversed)

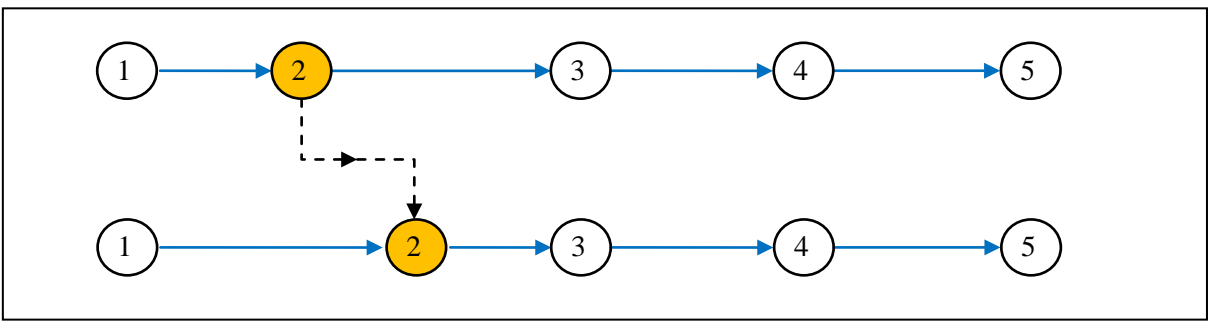

Figure 13.4 : Date shift (the date 2 is shifted between date 1 and date 3)

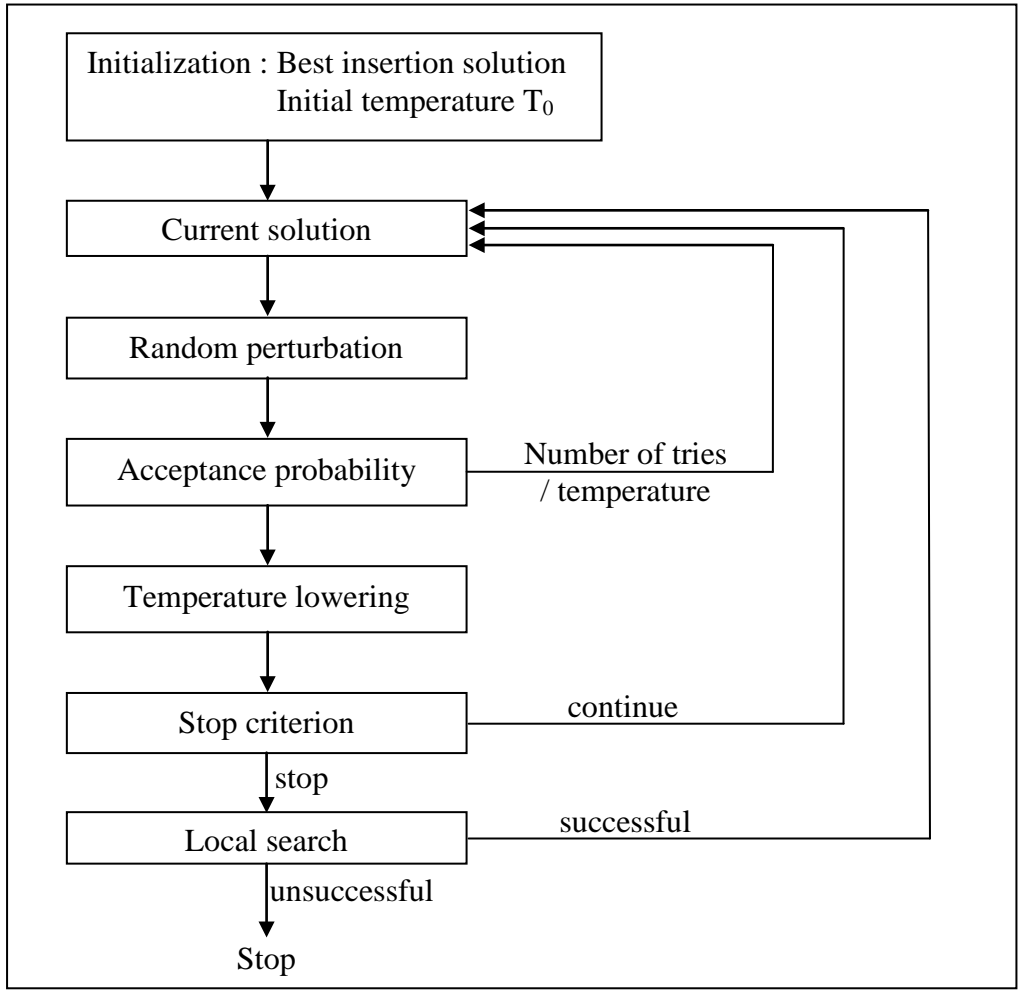

Figure 14 : Simulated annealing

The performance of the simulated annealing algorithm is checked on a sample of TSP instances used as tests benchmarks for stochastic programming[13] :
- Defi250 is a net contest with 250 fictitious towns
- Bier127 are the locations of 127 brasseries in München
- Lin318 is a laser pulsed drill with 318 holes
- Pcb442 is a Printed Circuit Board with 442 drills
- Att532 are the locations of the 532 US main towns

The Table 4 compares the results of the simulated annealing algorithm to the best published solutions[13]. The solutions are plotted on the Figure 15.

| Problem name | Cost found | Best known cost | Difference | Execution time |
|---|---|---|---|---|
| Defi250 | 11,9301 | 11,8092 | 1 % | 3 min |
| Bier127 | 118293 | 118282 | 0,009 % | 25 s |
| Lin318 | 42115 | 42029 | 0,2 % | 4 min |
| Pcb442 | 50927 | 50778 | 0,3 % | 11 min |
| Att532 | 28022 | 27686 | 1,2 % | 25 min |

Table 4 : TSP test cases

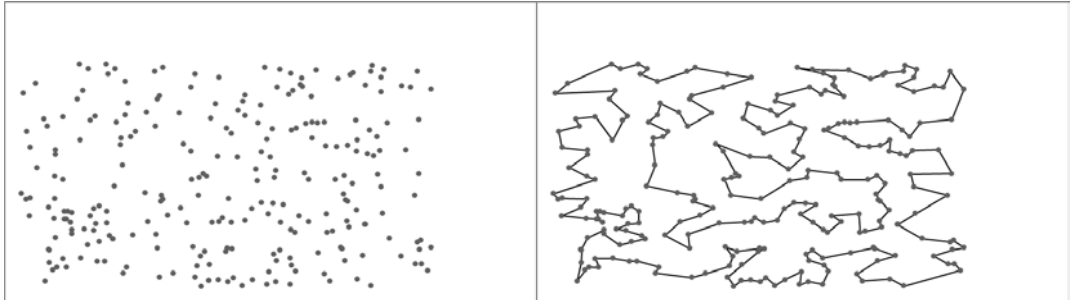

Figure 15.1 : Problem Defi250

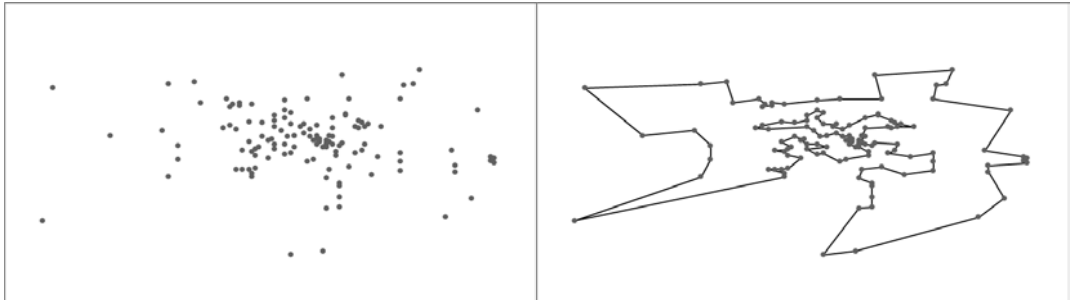

Figure 15.2 : Problem Bier127

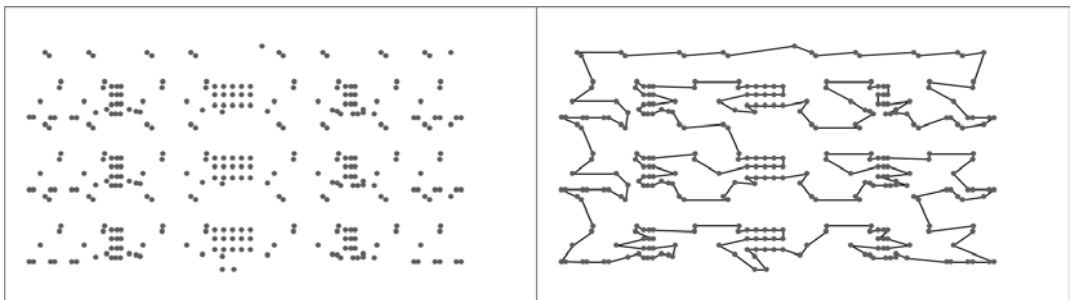

Figure 15.3 : Problem Lin318

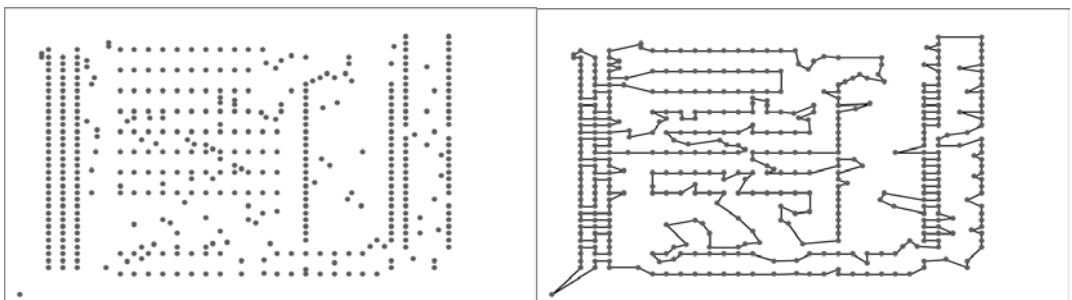

Figure 15.4 : Problem Pcb442

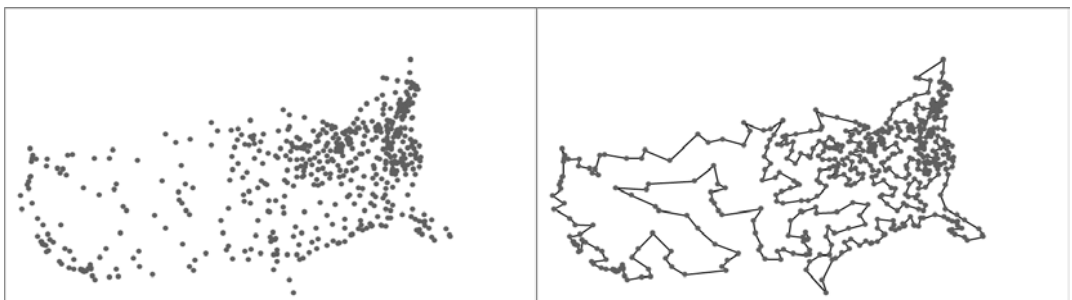

Figure 15.5 : Problem Att532

## 3.2 Response surface modelling

For the TSP, the cost function is simply the length of the closed path passing through the N nodes. The SDC cost function is more complex to assess :

- The edge valuations representing the transfers between debris are time-dependent.
- The nodes are gathered in sub-paths representing the successive missions.

The SDC cost function is assessed with the following procedure :

- A trial solution is still defined as a single path visiting the N candidates debris as for the TSP.
- A rendezvous date is associated to each node, with increasing values along the path.
- The cost of each edge is assessed by a response surface based on cost matrices.
- The path is sub-divided into m sub-paths of n nodes, representing the successive missions. The costs of these m sub-paths are denoted $K_i$, $1 \leq i \leq m$.
- The global cost function is the cost of the most expensive sub-path :
  $K = Max(K_1,…K_m)$

The cost evaluation procedure is recalled on the Figure 16 in the case of 3 missions of 5 debris to be selected in a list of 21 candidates.

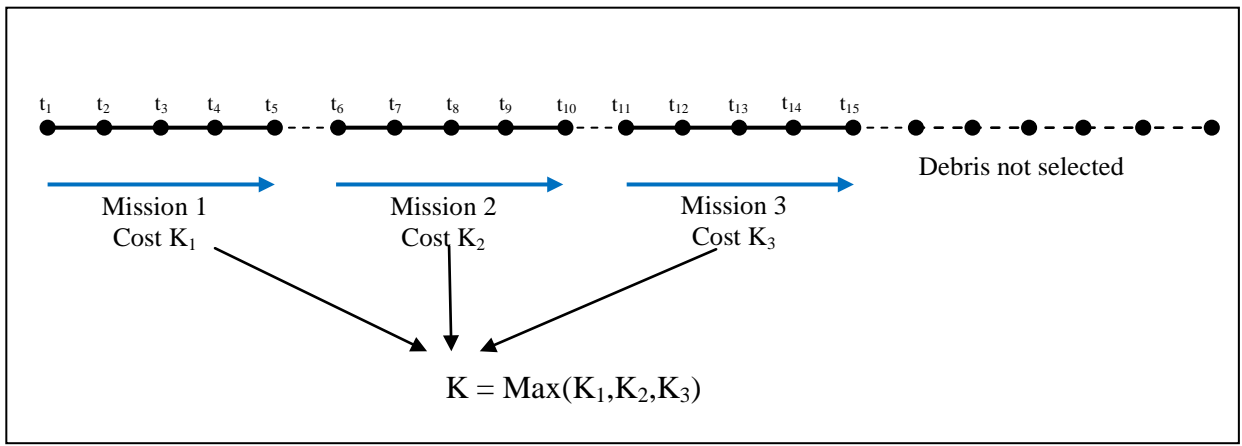

Figure 16 : SDC cost function

A trial solution is a path visiting the N debris. It is defined by :

- The successive debris numbers denoted $d_1, d_2, … , d_N$. The set $(d_1, d_2, … ,d_N)$ is a permutation of $(1,…,N)$.
- The successive rendezvous dates denoted $t_1, t_2, … , t_N$. These dates are increasing along the path : $t_1 < t_2 < … < t_{15}$ .

The main issue lies in the edge valuations which require solving a series of nonlinear transfer problems. An on-line optimization cannot be envisioned even with the simplified modelling presented in §2. It would result in huge computation times since millions of trial solutions are necessary before achieving the convergence of the simulated annealing algorithm.

In order to get a faster assessment the edges are valuated through a response surface modelling. The cost assessment is made by interpolation in the mesh of cost matrices spanning the cleaning program dates.

The pre-computed cost matrices $C(i,d,j,k)$ store the transfer costs from any debris j to any debris k≠j for a grid of discretized starting dates $\tau_i$ ($1 \leq i \leq n_t$) and discretized transfer durations $\Delta\tau_d$ ($1 \leq d \leq n_d$).
The $p^{th}$ edge on the trial path goes from the debris $d_p$ at the date $t_p$ to the debris $d_{p+1}$ at the date $t_{p+1}$. The transfer duration is denoted $\Delta t_p = t_{p+1} - t_p$. The interpolation consists in :

- Locating the starting date interval (index i) : $\tau_i \leq t_p < \tau_{i+1}$
- Locating the duration interval (index d) : $\Delta\tau_d \leq \Delta t_p < \Delta\tau_{d+1}$
- Selecting the matrices elements at the row $d_p$ (starting debris) and the column $d_{p+1}$ (arrival debris)
- Performing a bilinear interpolation on the intervals [$\tau_i$ ; $\tau_{i+1}$] and [$\Delta\tau_d$ ; $\Delta\tau_{d+1}$]

The interpolated cost of the $p^{th}$ edge is denoted $C_{int}(t_p, \Delta t_p, d_p, d_{p+1})$. The sub-path cost is the sum of the m edges interpolated costs. For the above illustration case (Figure 16), the 3 sub-paths representing the successive missions have the respective costs :

- Mission 1 : from the debris $d_1$ at the date $t_1$ to the debris $d_5$ at the date $t_5$

    $K_1 = C_{int}(t_1,\Delta t_1,d_1,d_2)\ \ \ \ + C_{int}(t_2,\Delta t_2,d_2,d_3)\ \ \ \ + C_{int}(t_3,\Delta t_3,d_3,d_4)\ \ \ \ + C_{int}(t_4,\Delta t_4,d_4,d_5)$

- Mission 2 : from the debris $d_6$ at the date $t_6$ to the debris $d_{10}$ at the date $t_{10}$

    $K_2 = C_{int}(t_6,\Delta t_6,d_6,d_7)\ \ \ \ + C_{int}(t_7,\Delta t_7,d_7,d_8)\ \ \ \ + C_{int}(t_8,\Delta t_8,d_8,d_9)\ \ \ \ + C_{int}(t_9,\Delta t_9,d_9,d_{10})$

- Mission 3 : from the debris $d_{11}$ at the date $t_{11}$ to the debris $d_{15}$ at the date $t_{15}$

    $K_3 = C_{int}(t_{11},\Delta t_{11},d_{11},d_{12}) + C_{int}(t_{12},\Delta t_{12},d_{12},d_{13}) + C_{int}(t_{13},\Delta t_{13},d_{13},d_{14}) + C_{int}(t_{14},\Delta t_{14},d_{14},d_{15})$

The total cost of the cleaning program is given by the most expensive mission :
$K = Max\ (K_1, K_2, K_3)$.

The cost assessment using this response surface modelling is sufficiently fast. It allows the application of the simulated annealing approach to the SDC problem. Some cautions are necessary in order to get reliable results :

- The mesh of discretized dates and durations must cover the cleaning program dates to avoid extrapolations (cf §2.3).
- The discretization step must be sufficiently small. The transfer cost function is indeed nonlinear with local minima and a too large time step may lead to skip good solutions.
- The cost found with the response surface modelling must be refined a posteriori by a simulation-based assessment. This is illustrated on the application case (§4).

## 3.3 Mission global cost

The actual cost function for the SDC problem is the fuel consumption per mission, which is the driver for the SDC vehicle design. The fuel consumption depends on the vehicle mass through the rocket equation. It is not an adequate measure for the transfer valuations since it depends on the transfer location along the path. Rather than the fuel consumption, the cost matrices store the velocity impulses which are intrinsic measures of the transfer costs.

The interpolated cost $C_{int}(t_p, \Delta t_p, d_p, d_{p+1})$ gives the required velocity impulse $\Delta V_p$ for the $p^{th}$ transfer going from the debris $d_p$ at the date $t_p$ to the debris $d_{p+1}$ at the date $t_{p+1} = t_p + \Delta t_p$.

The propellant consumed for this $p^{th}$ transfer is assessed from the rocket equation :

$$m_c(t_p) = M(t_p)\left(1 - e^{-\frac{\Delta V_p}{v_e}}\right)$$

$M(t_p)$ is the vehicle gross mass at the transfer beginning, $v_e$ is the engine exhaust velocity.

In addition to the transfer maneuvers, the mission assessment must also account for the debris operations, both in terms of durations and of released masses.
A fixed duration is allocated to each debris operations (observation, capture and deorbitation). This duration denoted $\Delta t_{oper}$ must be reserved within the duration of the cleaning program. For that purpose, it is directly taken into account when building the cost matrices by including a last waiting sequence of duration $\Delta t_{oper}$ in the transfer modelling, once the targeted debris is reached. An elementary transfer going from a debris 1 at the date $t_1$ to a debris 2 at the date $t_2$ is by this way completed at the date $t_2 - \Delta t_{oper}$. The operation durations between the successive transfers are thus implicitly accounted in the path valuation through the RSM.

The operation costs denoted $C_{oper}$ are stored on the cost matrix diagonals, so that they can be accounted in the mission global assessment. Two deorbitation options are envisioned :

- Either a deorbitation of the debris by the SDC vehicle
- Or an autonomous deorbitation of the debris using a "kit" supplied by the SDC vehicle.

The storage depends on the deorbitation option as follows.

- The first option consists in a deorbitation of the debris by the SDC vehicle. For each debris the velocity impulse $\Delta V_{oper}$ required by the deorbitation depends on the debris altitude and it can be assessed a priori. The deorbitation velocity impulses of the N debris are stored on the cost matrix diagonals. For an edge going from the debris j to the debris k, the deorbitation cost $\Delta V_{oper,k}$ of the debris k is added to the transfer interpolated cost, resulting in an additional fuel consumption.
- The second option consists in an autonomous deorbitation of the debris using a "kit" supplied by the SDC vehicle. The kit of mass $m_{oper}$ is attached to the debris, then the debris is released to perform the deorbitation maneuver. In that option, the masses of the N kits designed respectively for the N debris are stored on the cost matrix diagonals. At the end of the propelled transfer arriving on the debris k, the vehicle gross mass is decreased from the kit mass $m_{oper,k}$ released to the debris k.

The cost of the debris operations are thus taken into account within the path optimization, whatever the deorbitation option selected.

It is also possible to consider weights $w_k$ on the debris list to account for their priority, depending for example on their dangerousness. These weights come as multipliers on the cost matrix columns (arrival debris). The cost of an edge going from the debris j to the debris k is then assessed as :

$$w_k \left( C_{int}(t_p, \Delta t_p, j, k) + C_{oper}(k) \right)$$

### 3.4 Optimization process

The overall optimization process is split into three successive stages as presented on the Figure 6 :

- Cost matrices generation
- Path optimization
- Refined solution

The practical implementation of these three stages is described here after.

#### 3.4.1 Cost matrices generation

The debris orbits are retrieved at a given date from a database like the TLE of the NORAD[14]. A mesh of $n_t$ discretized starting dates and $n_d$ transfer durations is chosen in order to span the cleaning program forecast dates $[T_0 ; T_0 + \Delta T]$. The grid step results from a compromise between the response surface accuracy and the total computation time. The following choices are based on the mean duration per mission and on the mean duration per transfer. They have given an adequate compromise on the practical applications :

- $n_t \approx n \times m$ to associate one starting date per selected debris
  with  $\tau_1 = T_0$           (cleaning program starting date)
         $\tau_{nt} = T_0 + \Delta T$     (cleaning program ending date)

- $n_d \approx n$ to associate one transfer duration per sub-path debris
  with  $\Delta\tau_1 = \Delta T/m/2/n$   (minimum = half of mean transfer duration)
         $\Delta\tau_{nd} = \Delta T/m/2$    (maximum= half of the mean mission duration)

The $n_t \times n_d \times N \times (N-1)$ elementary transfer optimizations are run to fill the cost matrices.
These optimizations are independent from each other and they are parallelized on several processors. Each optimization is a NLP problem with 2 variables (drift orbit) and 1 constraint (final RAAN value). To spare some computation time, a filter discards the useless cases (ending date exceeding the end of the cleaning program) and the unfeasible transfers (requiring a drift altitude out of the allowed bounds). For such cases an arbitrary large cost value is stored in the corresponding matrix element so that it will not be selected on the path.

The optimization variables (drift orbit radius and inclination) are initialized automatically, depending on the debris relative RAAN values. The convergence is typically achieved in about 10 seconds.

### 3.4.2 Path optimization

The simulated annealing algorithm is applied, with the cost function assessed by a response surface modelling based on the cost matrices. The variables are the debris order and the rendezvous dates. The debris operations are accounted in terms of duration (taken into account in the cost matrices) and cost (stored on the matrix diagonals depending on the deorbitation option).

At the convergence, a solution path is issued defining the m missions of n debris each, with the optimized rendezvous dates. The convergence is achieved after some million trial solutions, with a typical computation time of a few minutes similarly to a classical TSP problem.

### 3.4.3 Refined solution

The cost function for the simulated annealing has been computed through a response surface modelling based on interpolations. The real cost is actually nonlinear and it can be significantly different of the RSM cost, depending on the mesh discretization. In order to check and refine the mission planning, and to get a reliable cost assessment, the missions are re-optimized using a simulation based software. The debris order is fixed, as well as the mission initial and final dates. For each mission, the rendezvous dates and the intermediate drift orbits are re-optimized to minimize the total $\Delta V$. The refined assessment yields the requirement for the SDC vehicle design.

## 4. Application case

The optimization method is illustrated on an application case with 21 debris, The cleaning program consists of 3 missions visiting 5 debris each one over a total duration of 45 months. The SDC vehicle shall use either a high thrust or a low thrust propulsion system.

### 4.1 Debris list

A list of 21 debris on circular orbits is considered, with the altitude ranging from 700 to 900 km, the inclination ranging from 97 to 99 deg, and the initial RAAN between 0 to 360 deg. The values of altitude, inclination and initial RAAN are uniformly distributed in their respective intervals. The Table 5 provides the orbital parameters of the 21 debris, with their nodal precession rate in the last column. The orbits are nearly sun-synchronous with precession rates close to the Sun precession rate (360 deg / 365.25 day = 0.986 deg/day).

For a real application case, the orbital parameters can be retrieved from official databases like the TLE of the NORAD.

| Debris number | Altitude (km) | Inclination (deg) | Initial RAAN (deg) | Precession rate (deg/day) |
|---|---|---|---|---|
| Debris 1 | 700 | 97.0 | 0. | 0,8429 |
| Debris 2 | 710 | 97.3 | 90. | 0,8745 |
| Debris 3 | 720 | 97.6 | 180. | 0,9058 |
| Debris 4 | 730 | 97.9 | 270. | 0,9367 |
| Debris 5 | 740 | 98.2 | 018. | 0,9672 |
| Debris 6 | 750 | 98.5 | 108. | 0,9975 |
| Debris 7 | 760 | 98.8 | 198. | 1,0273 |
| Debris 8 | 770 | 97.1 | 288. | 0,8260 |
| Debris 9 | 780 | 97.4 | 36. | 0,8565 |
| Debris 10 | 790 | 97.7 | 126. | 0,8866 |
| Debris 11 | 800 | 98.0 | 216. | 0,9165 |
| Debris 12 | 810 | 98.3 | 306. | 0,9460 |
| Debris 13 | 820 | 98.6 | 54. | 0,9752 |
| Debris 14 | 830 | 98.9 | 144. | 1,0040 |
| Debris 15 | 840 | 97.2 | 234. | 0,8094 |
| Debris 16 | 850 | 97.5 | 324. | 0,8389 |
| Debris 17 | 860 | 97.8 | 72. | 0,8681 |
| Debris 18 | 870 | 98.1 | 162. | 0,8969 |
| Debris 19 | 880 | 98.4 | 252. | 0,9254 |
| Debris 20 | 890 | 98.7 | 342. | 0,9536 |
| Debris 21 | 900 | 99.0 | 360. | 0,9815 |

Table 5 : List of 21 candidate debris

## 4.2 Cleaning program specification

The goal is to design the lightest vehicle able to perform the 3 successive missions. Each mission has to visit 5 debris, so that 15 debris out of the 21 candidates will be visited.

An average duration of 3 months per debris is considered, leading to a total duration of 45 months (1370 days) for the overall cleaning program. A 5 days duration is also allocated to each debris operations.

The altitude of the drift orbits is bounded between 400 km and 2000 km.

The mission consists in visiting successively the debris, without deorbitation maneuvers. The cost function is the total velocity impulse $\Delta V$ required for the orbital transfer maneuvers. The vehicle uses either a high thrust or a low thrust propulsion system. In the low thrust case an average acceleration level of 0.0035 m/s² is considered.

## 4.3 Cost matrices

The first stage of the solving method consists in building the cost matrices. Each cost matrix contains the transfer costs from any debris to any other for a given starting date and a given transfer duration. Filling one cost matrix requires solving 20×21 elementary transfer problems.

The matrices are assessed for a grid of discretized starting dates and durations, in order to span the total duration of the 3 missions. The application case specifies that 15 debris have to be visited within a 45 months period. We choose the following discretization :

- 16 starting dates ranging from 0 to 45 months
- 6 transfer durations ranging from 20 to 200 days

The total number of optimizations is 16×6×20×21 = 40320. Each optimization is achieved in about 10 seconds, leading to a total computation time of 112 hours. With a parallelization on 10 processors, the task is completed in a half day. The computation times are similar for the high thrust case and the low thrust case.

The date and duration grids and the cost matrices are written in an output file. The file is directly usable by the simulated annealing algorithm in order to build the response surface modelling.

## 4.4 Path optimization

The second stage of the solving method consists in finding the optimal mission planning leading to the minimal $\Delta V$ requirement per mission. The cost is assessed by interpolation in the cost matrices. The convergence of the simulated annealing is obtained in about 10 minutes with 200 million trials.

### 4.4.1 High thrust case

The Table 6 presents the cleaning program found with the simulated annealing algorithm, for the high thrust case. The corresponding mission planning is detailed in the Table 7.

It can be observed than the respective missions have close $\Delta V$ ranging from 820 m/s to 838 m/s. This balanced cost between the missions is an indication of a good behavior of the simulated annealing algorithm. The Figure 17 shows the convergence (cost function value, acceptation rate) plotted in function of the logarithm of the temperature.

|  | Dates (days) | Debris visited | Total ΔV (m/s) |
|---|---|---|---|
| Mission 1 | 0 - 545,3 | 16 - 20 - 21 - 5 - 17 | **820,0** |
| Mission 2 | 552,7 - 935,7 | 15 - 3 - 14 - 11 - 8 | **838,0** |
| Mission 3 | 942,1 - 1365,9 | 1 - 4 - 9 - 7 - 12 | **837,5** |

Table 6 : Cleaning program with high thrust (SA solution)

| Mission 1 | | | Mission 2 | | | Mission 3 | | |
|---|---|---|---|---|---|---|---|---|
| Debris number | Date (days) | ΔV (m/s) | Debris number | Date (days) | ΔV (m/s) | Debris number | Date (days) | ΔV (m/s) |
| **16** | 3,1 | 285,8 | **15** | 552,7 | 145,0 | **1** | 942,1 | 154,7 |
| **20** | 184,8 | 224,4 | **3** | 616,0 | 375,5 | **4** | 1014,6 | 362,8 |
| **21** | 375,0 | 216,0 | **14** | 771,5 | 140,4 | **9** | 1179,8 | 246,9 |
| **5** | 488,7 | 93,9 | **11** | 831,9 | 177,1 | **7** | 1268,0 | 73,2 |
| **17** | 545,3 | **820,0** | **8** | 935,7 | **838,0** | **12** | 1365,9 | **837,5** |

Table 7 : Mission planning with high thrust (SA solution)

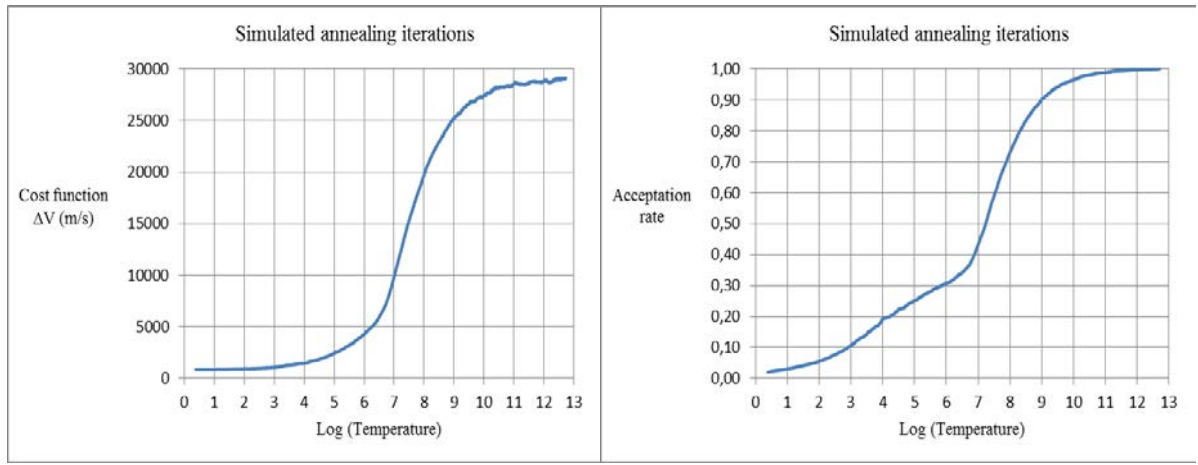

Figure 17 : Simulated annealing convergence

A sample of the iterations is listed in the Table 8, showing the path evolution. The last iterations show that several paths exist yielding close cost values.

| ΔV (m/s) | Debris order | | | | | | | | | | | | | | | | | | | | |
|---|---|---|---|---|---|---|---|---|---|---|---|---|---|---|---|---|---|---|---|---|---|
| 28353,3 | 15 | 14 | 13 | 12 | 4 | 10 | 9 | 8 | 7 | 6 | 5 | 11 | 3 | 1 | 2 | 16 | 17 | 18 | 19 | 20 | 21 |
| 22706,3 | 13 | 12 | 4 | 14 | 15 | 6 | 7 | 9 | 19 | 10 | 11 | 3 | 16 | 2 | 1 | 17 | 18 | 20 | 8 | 5 | 21 |
| 20859,4 | 13 | 12 | 7 | 14 | 15 | 4 | 9 | 10 | 19 | 17 | 11 | 3 | 16 | 2 | 1 | 18 | 20 | 8 | 5 | 6 | 21 |
| 20810,4 | 1 | 16 | 5 | 15 | 18 | 17 | 19 | 8 | 9 | 10 | 3 | 12 | 4 | 7 | 20 | 14 | 6 | 11 | 21 | 2 | 13 |
| 20644,6 | 1 | 16 | 8 | 15 | 18 | 9 | 5 | 19 | 17 | 10 | 3 | 12 | 4 | 7 | 20 | 14 | 6 | 11 | 21 | 2 | 13 |
| 18544,8 | 3 | 8 | 20 | 19 | 14 | 5 | 2 | 4 | 11 | 6 | 16 | 1 | 7 | 15 | 13 | 21 | 10 | 9 | 12 | 17 | 18 |
| 16611,0 | 3 | 7 | 20 | 19 | 14 | 5 | 8 | 4 | 1 | 6 | 11 | 16 | 2 | 15 | 13 | 21 | 10 | 9 | 12 | 17 | 18 |
| 16459,2 | 3 | 7 | 20 | 19 | 16 | 5 | 8 | 4 | 1 | 6 | 11 | 14 | 2 | 15 | 13 | 21 | 10 | 9 | 12 | 17 | 18 |
| 14129,7 | 18 | 7 | 20 | 8 | 1 | 12 | 2 | 6 | 10 | 16 | 19 | 14 | 4 | 15 | 13 | 21 | 9 | 5 | 17 | 3 | 11 |
| 13434,5 | 18 | 5 | 1 | 13 | 7 | 12 | 4 | 21 | 20 | 9 | 2 | 15 | 10 | 14 | 19 | 6 | 17 | 3 | 16 | 8 | 11 |
| | | | | | | | | | | | | | | | | | | | | | |
| 2016,4 | 19 | 11 | 7 | 16 | 4 | 15 | 6 | 3 | 18 | 8 | 20 | 2 | 21 | 5 | 10 | 12 | 1 | 9 | 13 | 17 | 14 |
| 1978,0 | 9 | 16 | 12 | 1 | 20 | 6 | 18 | 3 | 11 | 14 | 13 | 15 | 21 | 10 | 5 | 17 | 8 | 19 | 7 | 4 | 2 |
| 1817,0 | 13 | 18 | 14 | 3 | 15 | 16 | 4 | 11 | 8 | 6 | 9 | 12 | 17 | 2 | 20 | 19 | 7 | 10 | 21 | 5 | 1 |
| 1788,5 | 13 | 18 | 14 | 3 | 15 | 4 | 16 | 11 | 8 | 6 | 9 | 12 | 17 | 2 | 20 | 19 | 7 | 10 | 21 | 5 | 1 |
| 1769,6 | 16 | 12 | 1 | 20 | 9 | 4 | 7 | 19 | 11 | 8 | 13 | 18 | 15 | 5 | 10 | 14 | 6 | 2 | 17 | 21 | 3 |
| 1685,1 | 16 | 12 | 1 | 20 | 9 | 4 | 7 | 19 | 11 | 8 | 13 | 18 | 15 | 5 | 10 | 14 | 6 | 2 | 17 | 21 | 3 |
| 1604,2 | 10 | 14 | 18 | 6 | 15 | 4 | 16 | 7 | 8 | 11 | 2 | 5 | 21 | 20 | 17 | 12 | 19 | 13 | 1 | 9 | 3 |
| 1553,8 | 20 | 21 | 9 | 17 | 5 | 3 | 15 | 14 | 8 | 11 | 4 | 1 | 19 | 7 | 12 | 16 | 2 | 13 | 6 | 18 | 10 |
| 1514,0 | 10 | 6 | 17 | 5 | 21 | 16 | 4 | 1 | 7 | 12 | 15 | 13 | 18 | 3 | 8 | 2 | 19 | 14 | 9 | 11 | 20 |
| 1470,4 | 16 | 8 | 4 | 19 | 7 | 14 | 15 | 18 | 3 | 6 | 20 | 17 | 2 | 5 | 21 | 11 | 9 | 13 | 10 | 12 | 1 |
| 1374,6 | 12 | 16 | 8 | 19 | 7 | 18 | 6 | 3 | 11 | 14 | 17 | 20 | 2 | 5 | 21 | 15 | 1 | 10 | 9 | 4 | 13 |
| 1326,3 | 2 | 6 | 18 | 3 | 14 | 16 | 19 | 7 | 1 | 4 | 13 | 15 | 21 | 5 | 10 | 8 | 9 | 20 | 17 | 11 | 12 |
| 1301,8 | 2 | 6 | 18 | 3 | 14 | 19 | 16 | 7 | 1 | 4 | 13 | 15 | 21 | 10 | 5 | 8 | 9 | 20 | 17 | 11 | 12 |
| 1268,6 | 19 | 8 | 4 | 12 | 1 | 15 | 3 | 6 | 11 | 14 | 20 | 17 | 21 | 5 | 10 | 9 | 2 | 13 | 18 | 16 | 7 |
| 1251,0 | 3 | 7 | 19 | 4 | 16 | 15 | 14 | 11 | 8 | 6 | 17 | 20 | 21 | 5 | 10 | 9 | 12 | 13 | 18 | 2 | 1 |
| 1239,8 | 11 | 15 | 3 | 18 | 6 | 4 | 16 | 7 | 12 | 9 | 17 | 20 | 21 | 10 | 5 | 13 | 2 | 8 | 19 | 14 | 1 |
| 1231,6 | 10 | 6 | 14 | 3 | 11 | 4 | 16 | 19 | 9 | 12 | 18 | 13 | 15 | 5 | 21 | 17 | 20 | 8 | 7 | 2 | 1 |
| 1213,0 | 15 | 3 | 14 | 18 | 6 | 19 | 7 | 16 | 4 | 1 | 20 | 17 | 21 | 5 | 10 | 9 | 13 | 11 | 8 | 2 | 12 |
| 1210,4 | 15 | 3 | 14 | 18 | 6 | 19 | 7 | 16 | 4 | 1 | 17 | 20 | 21 | 5 | 10 | 9 | 13 | 11 | 8 | 2 | 12 |
| 1142,3 | 12 | 16 | 4 | 19 | 7 | 3 | 14 | 11 | 8 | 6 | 13 | 21 | 5 | 10 | 15 | 9 | 18 | 17 | 1 | 20 | 2 |
| 1128,2 | 17 | 13 | 6 | 18 | 3 | 16 | 4 | 7 | 12 | 9 | 2 | 21 | 5 | 10 | 15 | 19 | 14 | 20 | 11 | 1 | 8 |
| 1119,5 | 7 | 11 | 15 | 3 | 14 | 12 | 1 | 4 | 16 | 19 | 20 | 2 | 21 | 5 | 10 | 13 | 9 | 17 | 6 | 8 | 18 |
| 1105,9 | 12 | 16 | 4 | 7 | 19 | 18 | 14 | 11 | 8 | 6 | 13 | 15 | 10 | 21 | 5 | 9 | 2 | 1 | 17 | 3 | 20 |
| 1052,4 | 12 | 16 | 4 | 7 | 19 | 3 | 14 | 11 | 8 | 6 | 13 | 15 | 10 | 21 | 5 | 9 | 2 | 1 | 17 | 18 | 20 |
| 996,7 | 12 | 16 | 4 | 19 | 7 | 3 | 14 | 11 | 8 | 6 | 13 | 15 | 21 | 5 | 10 | 9 | 2 | 1 | 17 | 18 | 20 |
| 962,6 | 12 | 8 | 4 | 19 | 7 | 18 | 6 | 3 | 11 | 14 | 13 | 15 | 10 | 5 | 21 | 1 | 2 | 20 | 9 | 16 | 17 |
| 958,0 | 12 | 4 | 8 | 19 | 7 | 18 | 6 | 3 | 11 | 14 | 13 | 15 | 10 | 5 | 21 | 1 | 2 | 20 | 9 | 16 | 17 |
| 946,7 | 18 | 3 | 14 | 15 | 6 | 16 | 7 | 4 | 12 | 9 | 17 | 20 | 21 | 5 | 10 | 8 | 11 | 19 | 13 | 2 | 1 |
| 900,6 | 6 | 18 | 14 | 3 | 15 | 19 | 16 | 7 | 4 | 1 | 17 | 20 | 21 | 5 | 10 | 13 | 2 | 8 | 12 | 11 | 9 |
| 897,8 | 6 | 18 | 14 | 3 | 15 | 19 | 16 | 7 | 4 | 1 | 17 | 20 | 21 | 5 | 10 | 13 | 9 | 8 | 2 | 11 | 12 |
| 891,2 | 6 | 18 | 14 | 3 | 15 | 19 | 16 | 7 | 4 | 1 | 20 | 17 | 21 | 5 | 10 | 13 | 9 | 8 | 11 | 2 | 12 |
| 889,2 | 2 | 13 | 17 | 5 | 21 | 15 | 3 | 14 | 11 | 8 | 1 | 4 | 9 | 7 | 12 | 19 | 16 | 18 | 6 | 10 | 20 |
| 884,0 | 2 | 13 | 17 | 5 | 21 | 15 | 3 | 14 | 11 | 8 | 1 | 4 | 9 | 7 | 12 | 19 | 16 | 18 | 6 | 10 | 20 |
| 881,7 | 2 | 13 | 17 | 5 | 21 | 15 | 3 | 14 | 11 | 8 | 1 | 4 | 9 | 7 | 12 | 19 | 16 | 18 | 6 | 10 | 20 |
| 879,3 | 1 | 21 | 5 | 9 | 20 | 15 | 3 | 14 | 11 | 8 | 16 | 19 | 4 | 7 | 12 | 18 | 6 | 13 | 2 | 17 | 10 |
| 871,8 | 1 | 21 | 5 | 9 | 20 | 15 | 3 | 14 | 11 | 8 | 16 | 19 | 4 | 7 | 12 | 18 | 6 | 13 | 2 | 17 | 10 |
| **838,0** | 16 | 20 | 21 | 5 | 17 | 15 | 3 | 14 | 11 | 8 | 1 | 4 | 9 | 7 | 12 | 19 | 2 | 13 | 18 | 6 | 10 |

Table 8 : Simulated annealing iterations (sample)

### 4.4.2 Low thrust case

The Table 9 presents the cleaning program found with the simulated annealing algorithm, for the low thrust case. The corresponding mission planning is detailed in the Table 10.
As for the high thrust case a good balance can be observed between the successive missions with the Δvelocity impulses ranging from 963 m/s to 973 m/s.

The orbital maneuvers using a low thrust engine take a significant time that can no longer be used for the drift phase. The price for the shortened drift duration is paid by a "farther" drift orbit in order to accelerate the RAAN precession, so that the RAAN constraint can still be met in the prescribed transfer duration. Also the low thrust transfer incurs velocity losses that do not exist in the impulsive modelling.
The mission costs with a low thrust engine are therefore slightly higher in terms of velocity impulse, but significant savings are expected in terms of fuel consumption due to a much better exhaust velocity.

|  | Dates (days) | Debris visited | Total ΔV (m/s) |
|---|---|---|---|
| Mission 1 | 0 - 556,9 | 6 - 14 - 18 - 3 - 15 | **969,2** |
| Mission 2 | 563,8 - 978,7 | 19 - 7 - 16 - 4 - 1 | **963,7** |
| Mission 3 | 987,6 - 1361,0 | 17 - 20 - 21 - 5 - 10 | **972,4** |

Table 9 : Cleaning program with low thrust (SA solution)

| Mission 1 | | | Mission 2 | | | Mission 3 | | |
|---|---|---|---|---|---|---|---|---|
| Debris number | Date (days) | ΔV (m/s) | Debris number | Date (days) | ΔV (m/s) | Debris number | Date (days) | ΔV (m/s) |
| **6** | 0,7 | 377,7 | **19** | 563,8 | 240,7 | **17** | 987,6 | 220,2 |
| **14** | 174,1 | 184,9 | **7** | 621,0 | 276,6 | **20** | 1047,4 | 443,5 |
| **18** | 302,9 | 262,3 | **16** | 714,1 | 274,7 | **21** | 1252,3 | 188,1 |
| **3** | 483,8 | 144,4 | **4** | 886,3 | 171,6 | **5** | 1285,8 | 120,6 |
| **15** | 556,9 | **969,2** | **1** | 978,7 | **963,7** | **10** | 1361,0 | **972,4** |

Table 10 : Mission planning with low thrust (SA solution)

Comparing with the high thrust case, it can be observed that some debris remain gathered, for example (5 – 17 – 20 – 21), (3 – 14 – 15) and (1 – 4 – 7). A part of the explanation lies in the close initial RAAN values (Table 13 and Table 16)

### 4.5 Refined solution

In order to get a reliable cost assessment, the three missions are re-optimized using a simulation based software. The debris order is fixed, as well as the mission initial and final dates. For each mission, the rendezvous dates and the intermediate drift orbits are re-optimized to minimize the total mission ΔV. The rendezvous with the successive debris is constrained with a maximum RAAN deviation of 1 degree.

### 4.5.1 High thrust case

The Table 11 presents the cleaning program found after the dates and maneuvers re-optimization, for the high thrust case. The total ΔV per mission are given in the last column (in parenthesis the previous RSM assessment issued from the simulated annealing).
The corresponding mission planning is detailed in the Table 12.

|  | Dates (days) | Debris visited | Total ΔV (m/s) |
|---|---|---|---|
| Mission 1 | 0 - 545,3 | 16 - 20 - 21 - 5 - 17 | **811.1** (820.0) |
| Mission 2 | 552,7 - 935,7 | 15 - 3 - 14 - 11 - 8 | **711.9** (838.0) |
| Mission 3 | 942,1 - 1365,9 | 1 - 4 - 9 - 7 - 12 | **785.1** (837.5) |

Table 11 : Cleaning program with high thrust (simulation)

| Mission 1 | | | Mission 2 | | | Mission 3 | | |
|---|---|---|---|---|---|---|---|---|
| Debris number | Date (days) | ΔV (m/s) | Debris number | Date (days) | ΔV (m/s) | Debris number | Date (days) | ΔV (m/s) |
| **16** | 3,1 | 287,1 | **15** | 552,7 | 141,5 | **1** | 942,1 | 119,2 |
| **20** | 183,1 | 210,8 | **3** | 563,3 | 291,8 | **4** | 976,8 | 411,8 |
| **21** | 389,3 | 202,2 | **14** | 781,7 | 132,2 | **9** | 1143,4 | 183,5 |
| **5** | 513,6 | 111,0 | **11** | 823,0 | 146,4 | **7** | 1177,3 | 70,6 |
| **17** | 545,4 | **811,1** | **8** | 935,8 | **711,9** | **12** | 1365,8 | **785,1** |

Table 12 : Mission planning with high thrust (simulation)

The comparison to the simulated annealing results using the RSM cost assessment shows an improvement on all the missions owing to the drift orbit parameters refined optimization. The main changes are marked by the orange colored cells.

- The improvement on the mission 2 comes mainly from the second leg whose ΔV is reduced of 85 m/s. This is explained by the starting date advance (-50 days) which increases the transfer duration up to 218 days. Such a solution could not be detected with the response surface modelling, because the mesh was defined with a maximum transfer duration of 200 days. In such a case, it could be useful to iterate the whole process with an updated discretization.

- The improvement on the mission 3 comes mainly from the advance of the second rendezvous date (-40 days). The ΔV of the second leg is increased of about 50 m/s, but this is globally counterbalanced by the gains on the first and the third legs. A refined mesh discretization may help capturing these nonlinearities of the cost function within the RSM.

Another interesting observation can be done by considering for each mission the RAAN of the selected debris at the mission starting date (Table 13). It can be noted that the debris assigned to each mission are somewhat gathered by their initial RAAN values, in increasing order. This is not an absolute rule, since the optimal order also depends on the other orbital parameters (radius and inclination). When large RAAN differences exist at the mission starting date (for example for the third and fifth debris of the mission 3), they are nullified by long transfer durations (up to 6 months).

| Mission 1 ($t_1$ = 3,1 days) | | Mission 2 ($t_6$ = 552,7 days) | | Mission 3 ($t_{11}$ = 942,1 days) | |
|---|---|---|---|---|---|
| Debris number | RAAN (deg) | Debris number | RAAN (deg) | Debris number | RAAN (deg) |
| **16** | -33,4 | **15** | -38,6 | **1** | 74,1 |
| **20** | -15,1 | **3** | -39,4 | **4** | 72,4 |
| **21** | 3,0 | **14** | -21,1 | **9** | 122,8 |
| **5** | 21,0 | **11** | 2,5 | **7** | 74,1 |
| **17** | 74,7 | **8** | 24,5 | **12** | 117,2 |

Table 13 : RAAN values at the mission beginning (high thrust)

The re-optimized mission planning is detailed in the Table 14. The last column checks the RAAN constraint at the successive rendezvous with the debris. Some observations can be done on the green and orange colored cells.

- The short drift durations (green cells) correspond to a drift orbit close to the starting orbit. The natural precession on the initial orbit is nearly sufficient to reach the targeted RAAN value, so that there is no need for significantly changing the precession rate (cf §2.2.3). This reduces the ΔV required for the transfer.

- The large drift durations (orange cells) correspond to low altitude drift orbits. These transfers require a large RAAN change achieved by both an accelerated precession rate (low altitude) and a long duration (about 6 months). The costs of these legs represent about half the mission cost.

| Mission 1 | | | SDC vehicle | | Transfer | | Drift orbit | | Debris orbit | | RAAN (deg) | |
|---|---|---|---|---|---|---|---|---|---|---|---|---|
| Edge number | Initial debris | Final debris | Date (days) | Total DV (m/s) | Duration (days) | DV (m/s) | Altitude (km) | Inclination (deg) | Altitude (km) | Inclination (deg) | SDC Vehicle | Debris |
| 1 |  | 16 | 3,1 | 0 | 3,1 | 0 |  |  | 850,0 | 97,50 | 326,6 | 326,6 |
| 2 | 16 | 20 | 183,1 | 287,1 | 180,0 | 287,1 | 708,0 | 98,84 | 890,0 | 98,70 | 156,7 | 156,6 |
| 3 | 20 | 21 | 389,3 | 497,9 | 206,2 | 210,8 | 715,8 | 99,20 | 900,0 | 99,00 | 22,2 | 22,1 |
| 4 | 21 | 5 | 513,6 | 700,1 | 124,3 | 202,2 | 695,4 | 98,90 | 740,0 | 98,20 | 155,0 | 154,8 |
| 5 | 5 | 17 | 545,4 | 811,1 | 31,8 | 111,0 | 712,9 | 98,24 | 860,0 | 97,80 | 185,7 | 185,4 |

| Mission 2 | | | SDC vehicle | | Transfer | | Drift orbit | | Debris orbit | | RAAN (deg) | |
|---|---|---|---|---|---|---|---|---|---|---|---|---|
| Edge number | Initial debris | Final debris | Date (days) | Total DV (m/s) | Duration (days) | DV (m/s) | Altitude (km) | Inclination (deg) | Altitude (km) | Inclination (deg) | SDC Vehicle | Debris |
| 1 |  | 15 | 552,7 | 0 | 3,1 | 0 |  |  | 840,0 | 97,20 | 321,7 | 321,4 |
| 2 | 15 | 3 | 563,3 | 141,5 | 10,6 | 141,5 | 830,3 | 96,93 | 720,0 | 97,60 | 330,5 | 330,2 |
| 3 | 3 | 14 | 781,7 | 433,3 | 218,4 | 291,8 | 572,5 | 98,55 | 830,0 | 98,90 | 209,2 | 208,8 |
| 4 | 14 | 11 | 823,0 | 565,5 | 41,3 | 132,2 | 812,4 | 98,93 | 800,0 | 98,00 | 250,5 | 250,2 |
| 5 | 11 | 8 | 935,8 | 711,9 | 112,9 | 146,4 | 825,9 | 97,10 | 770,0 | 97,10 | 341,4 | 340,9 |

| Mission 3 | | | SDC vehicle | | Transfer | | Drift orbit | | Debris orbit | | RAAN (deg) | |
|---|---|---|---|---|---|---|---|---|---|---|---|---|
| Edge number | Initial debris | Final debris | Date (days) | Total DV (m/s) | Duration (days) | DV (m/s) | Altitude (km) | Inclination (deg) | Altitude (km) | Inclination (deg) | SDC Vehicle | Debris |
| 1 |  | 1 | 942,1 | 0 | 3,1 | 0 |  |  | 700,0 | 97,00 | 74,7 | 74,1 |
| 2 | 1 | 4 | 976,8 | 119,2 | 34,8 | 119,2 | 702,9 | 97,32 | 730,0 | 97,90 | 105,4 | 105,0 |
| 3 | 4 | 9 | 1143,4 | 531,0 | 166,6 | 411,8 | 412,4 | 98,29 | 780,0 | 97,40 | 295,9 | 295,3 |
| 4 | 9 | 7 | 1177,3 | 714,5 | 33,9 | 183,5 | 759,2 | 98,07 | 760,0 | 98,80 | 328,0 | 327,5 |
| 5 | 7 | 12 | 1365,8 | 785,1 | 188,6 | 70,6 | 763,7 | 98,69 | 810,0 | 98,30 | 158,8 | 158,1 |

Table 14 : Re-optimized missions with high thrust

### 4.5.2 Low thrust case

The Table 15 presents the cleaning program found after the dates and maneuvers re-optimization, for the low thrust case. The total ΔV per mission are given in the last column (in parenthesis the previous RSM assessment coming from the simulated annealing).
The corresponding mission planning is detailed in the Table 16.

|  | Dates (days) | Debris visited | Total ΔV (m/s) |
|---|---|---|---|
| Mission 1 | 0 - 558,1 | 6 - 14 - 18 - 3 - 15 | **903.6**   (969.2) |
| Mission 2 | 563,8 - 977,5 | 19 - 7 - 16 - 4 - 1 | **921.2**   (963.7) |
| Mission 3 | 987,6 - 1359,9 | 17 - 20 - 21 - 5 - 10 | **926.0**   (972.4) |

Table 15 : Cleaning program with low thrust (simulation)

| Mission 1 | | | Mission 2 | | | Mission 3 | | |
|---|---|---|---|---|---|---|---|---|
| Debris number | Date (days) | ΔV (m/s) | Debris number | Date (days) | ΔV (m/s) | Debris number | Date (days) | ΔV (m/s) |
| **6** | 0,7 | 369,8 | **19** | 563,8 | 266,1 | **17** | 987,6 | 197,8 |
| **14** | 178,8 | 194,3 | **7** | 622,3 | 247,4 | **20** | 1059,0 | 440,4 |
| **18** | 243,4 | 219,2 | **16** | 680,3 | 222,3 | **21** | 1262,9 | 182,3 |
| **3** | 486,6 | 120,3 | **4** | 954,3 | 185,4 | **5** | 1272,0 | 105,5 |
| **15** | 558,1 | **903,6** | **1** | 977,5 | **921,2** | **10** | 1359,9 | **926,0** |

Table 16 : Mission planning with low thrust (simulation)

Similarly to the high thrust case, the cost of all the missions are improved wrt the simulated annealing results using the RSM. The main changes are marked by the orange colored cells.

- The improvement on the mission 2 comes mainly from the third leg (-40 m/s) by advancing the transfer starting date (-60 days). The transfer duration is increased up to 243 days. This solution was not detected with the response surface modelling, because the mesh was defined with a maximum transfer duration of 200 days..

- The improvement on the mission 3 comes mainly from the delay of the fourth rendezvous date (+68 days). The ΔV of the third leg is reduced of about 50 m/s owing to a very long transfer duration (274 days).

As for the high thrust case, a gathering of the debris wrt to their RAAN values at the successive mission beginning can be observed (Table 17).

| Mission 1 ($t_1$ = 3,1 days) | | Mission 2 ($t_6$ = 563,8 days) | | Mission 3 ($t_{11}$ = 987,6 days) | |
|---|---|---|---|---|---|
| Debris number | RAAN (deg) | Debris number | RAAN (deg) | Debris number | RAAN (deg) |
| **6** | 108,7 | **19** | 53,7 | **17** | 209,3 |
| **14** | 144,7 | **7** | 57,2 | **20** | 203,7 |
| **18** | 162,6 | **16** | 76,9 | **21** | 249,3 |
| **3** | 180,6 | **4** | 78,1 | **5** | 253,2 |
| **15** | 234,6 | **1** | 115,2 | **10** | 281,6 |

Table 17 : RAAN values at the mission beginning (low thrust)

The re-optimized mission planning is detailed in the Table 18. The last column checks the RAAN constraint at the successive rendezvous with the debris. The same observations can be done as for the high thrust case, regarding the short drift durations (drift orbit close to the starting orbit, reachable with a low ΔV) and the large drift durations (low altitude drift orbits inducing a high ΔV).

With the mean acceleration level considered for this application case (0.0035 m/s²) is considered the propelled transfer durations do not exceed one day. The Edelbaum based transfer modelling, and the generic transfer strategy using a drift orbit are nearly optimal in this frame.

| Mission 1 | | | SDC vehicle | | Transfer | | Drift orbit | | Debris orbit | | RAAN (deg) | |
|---|---|---|---|---|---|---|---|---|---|---|---|---|
| Edge number | Initial debris | Final debris | Date (days) | Total DV (m/s) | Duration (days) | DV (m/s) | Altitude (km) | Inclination (deg) | Altitude (km) | Inclination (deg) | SDC Vehicle | Debris |
| 1 |   | 6 | 0,7 | 0 | 0,7 | 0 |   |   | 750,0 | 98,50 | 108,7 | 108,7 |
| 2 | 6 | 14 | 178,8 | 369,8 | 178,1 | 369,8 | 522,8 | 99,23 | 830,0 | 98,90 | 323,6 | 323,5 |
| 3 | 14 | 18 | 243,4 | 564,1 | 64,6 | 194,3 | 943,0 | 98,21 | 870,0 | 98,10 | 20,5 | 20,3 |
| 4 | 18 | 3 | 486,6 | 783,3 | 243,2 | 219,2 | 619,1 | 97,91 | 720,0 | 97,60 | 261,1 | 260,8 |
| 5 | 3 | 15 | 558,1 | 903,6 | 71,5 | 120,3 | 701,9 | 97,61 | 840,0 | 97,20 | 326,0 | 325,8 |

| Mission 2 | | | SDC vehicle | | Transfer | | Drift orbit | | Debris orbit | | RAAN (deg) | |
|---|---|---|---|---|---|---|---|---|---|---|---|---|
| Edge number | Initial debris | Final debris | Date (days) | Total DV (m/s) | Duration (days) | DV (m/s) | Altitude (km) | Inclination (deg) | Altitude (km) | Inclination (deg) | SDC Vehicle | Debris |
| 1 |   | 19 | 563,8 | 0 | 0,7 | 0 |   |   | 880,0 | 98,40 | 54,0 | 53,7 |
| 2 | 19 | 7 | 622,3 | 266,1 | 58,5 | 266,1 | 581,7 | 98,57 | 760,0 | 98,80 | 117,7 | 117,3 |
| 3 | 7 | 16 | 680,3 | 513,5 | 58,0 | 247,4 | 760,1 | 98,60 | 850,0 | 97,50 | 175,0 | 174,7 |
| 4 | 16 | 4 | 954,3 | 735,8 | 274,0 | 222,3 | 594,4 | 97,75 | 730,0 | 97,90 | 84,4 | 83,9 |
| 5 | 4 | 1 | 977,5 | 921,2 | 23,2 | 185,4 | 720,1 | 97,31 | 700,0 | 97,00 | 104,5 | 104,0 |

| Mission 3 | | | SDC vehicle | | Transfer | | Drift orbit | | Debris orbit | | RAAN (deg) | |
|---|---|---|---|---|---|---|---|---|---|---|---|---|
| Edge number | Initial debris | Final debris | Date (days) | Total DV (m/s) | Duration (days) | DV (m/s) | Altitude (km) | Inclination (deg) | Altitude (km) | Inclination (deg) | SDC Vehicle | Debris |
| 1 |   | 17 | 987,6 | 0 | 0,7 | 0 |   |   | 860,0 | 97,80 | 209,9 | 209,3 |
| 2 | 17 | 20 | 1059,0 | 197,8 | 71,4 | 197,8 | 945,4 | 98,14 | 890,0 | 98,70 | 272,4 | 271,8 |
| 3 | 20 | 21 | 1262,9 | 638,2 | 204,0 | 440,4 | 507,1 | 99,22 | 900,0 | 99,00 | 160,2 | 159,5 |
| 4 | 21 | 5 | 1272,0 | 820,5 | 9,1 | 182,3 | 749,7 | 98,27 | 740,0 | 98,20 | 169,0 | 168,3 |
| 5 | 5 | 10 | 1359,9 | 926,0 | 87,9 | 105,5 | 750,2 | 98,12 | 790,0 | 97,70 | 252,4 | 251,7 |

Table 18 : Re-optimized missions with low thrust

## 5. Conclusion

In order to clean the LEO region from the most dangerous debris, a cleaning program is envisioned. It consists of several successive missions, performed by similar vehicles, in order to achieve a mean removal rate of 5 debris per year.

A solution method is proposed for the planning of these successive Space Debris Collecting missions. The goal is to minimize the fuel required by the most expensive mission, with the perspective of designing a generic vehicle compliant of the cleaning program.

The problem mixes combinatorial optimization to select and order the debris among a list of candidates, and continuous optimization to fix the rendezvous dates and to define the minimum fuel orbital maneuvers. The solution method proposed consists in three stages.

Firstly the orbital transfer problem is simplified by considering a generic transfer strategy suited either to a high thrust or a low thrust vehicle. A response surface modelling is built by solving the reduced problem for all pairs of debris and for discretized dates, and storing the results in cost matrices. This first stage is parallelized on several processors. The results of this series of optimizations are stored in cost matrices.

Secondly a simulated annealing algorithm is applied to find the optimal mission planning. The cost function is assessed by interpolation on the response surface based on the cost matrices. This allows the convergence of the simulated algorithm in a limited computation time, yielding an optimal mission planning.

Thirdly the successive missions are re-optimized in terms of transfer maneuvers and dates without changing the debris order. This continuous control problem is simulation-based, taking into account the problem nonlinearities that were not captured by the response surface modelling. It yields a refined solution with the performance requirement for designing the future Space Debris Collecting vehicle.

The method is applicable for a large list of debris and for various assumptions regarding the cleaning program (number of missions, number of debris per mission, total duration, deorbitation scenario, high or low thrust vehicle). Weights are attributed to the debris in order to account for their dangerousness and to assign if desired a priority in the selection process. The generic transfer strategy can be considered as near optimal as long as a significant duration is allocated to the drift phases. For a low thrust vehicle, the available acceleration level must be sufficient, so that the propelled transfers do not exceed a few days.

The overall optimization process is automatized and a mission planning can be established in a few hours. It is exemplified on an application case with 3 missions to plan, each mission visiting 5 SSO debris to be selected in a list of 21 candidates. A low thrust propulsion leads to a slightly higher $\Delta V$ requirement, due partly to the reduced durations allocated to the drift phases, and partly to the velocity losses that do not exist in the impulsive modelling. Nevertheless significant savings are expected in terms of fuel consumption due to a much better exhaust velocity.

## Acknowledgment

This work was carried out at Airbus Defence and Space in 2013-2014 in the frame of the internal R&D. I would like to thank the R&D team project for having supported this work.

# Acronyms

| | |
|---|---|
| LEO | Low Earth Orbit |
| SSO | Sun Synchronous Orbit |
| RAAN | Right Ascension of Ascending Node |
| SDC | Space Debris Collecting |
| TSP | Travelling Salesman Problem |
| RSM | Response Surface Modelling |
| SA | Simulated Annealing |
| NLP | NonLinear Programming |

# Summary

## List of Tables

## List of Figures